\documentclass[notitlepage,leqno,10pt]{article}
\usepackage{amsmath}
\usepackage{amsfonts}
\usepackage{amssymb}
\title{ Hypersurfaces with free boundary and large  constant mean curvature: concentration along submanifolds}
\author{Mouhamed Moustapha  Fall$^{\rm a}$  and Fethi Mahmoudi$^{\rm a, b}$}
\date{}
\newtheorem{theorem}{Theorem}[section]
\newtheorem{proposition}{Proposition}[section]

\newtheorem{lemma}{Lemma}[section]

\newtheorem{remark}{Remark}[section]
\newcommand{\R}{\mathbb{R}}
\newcommand{\N}{\mathbb{N}}

\newcommand{\al}{\alpha}
\newcommand{\be}{\beta}
\newcommand{\e}{\varepsilon}
\newcommand{\del}{\partial}

\newcommand{\Om}{{\Omega}}
\newcommand{\calC}{{\mathcal C}}
\newcommand{\calE}{{\mathcal E}}

\newcommand{\calL}{\mathcal{ L}}
\newcommand{\calM}{{\mathcal M}}
\newcommand{\calO}{{\mathcal O}}

\newcommand{\calV}{{\mathcal V}}

\newcommand{\s}{\sigma}
\newcommand{\ov}{\overline}
\newcommand{\oy}{\overline{y}}
\newcommand{\z}{\zeta}
\newcommand{\n }{\nabla }

\newcommand{\pa}{{\partial}}

\newcommand{\de}{{\delta}}
\newcommand{\T}{{\Theta}}
\newcommand{\te}{{\theta}}
\newcommand{\TT}{{\tilde{\Theta}}}
\newcommand{\Tm}{{\Theta^{n+1}}}
\newcommand{\ra}{{\rangle}}
\newcommand{\la}{{\langle}}
\newcommand{\bfp}{\mathbf{p}}
\begin{document}

\maketitle

\begin{center}

$^{\rm a}${\small SISSA, Sector of Mathematical Analysis \\ Via
Beirut 2-4, 34014 Trieste, Italy}\\[2mm]
$^{\rm b}${\small D\'epartement de Mat\'ematiques, Facult\'e des sciences de Tunis\\
Campus Universitaire 2092 Tunis, Tunisia.  }
\end{center}

\footnotetext[1]{E-mail addresses:  fall@sissa.it (M.M. Fall), mahmoudi@ssissa.it
(F.Mahmoudi)}

\

\noindent {\sc abstract}. Given a domain $\Om$ of $\R^{m+1}$ and a $k$-dimensional
non-degenerate minimal submanifold $K$ of $\pa \Om$ with $1\le k\le m-1$, we prove
the existence of a family of embedded constant mean curvature hypersurfaces which as
their mean curvature tends to infinity concentrate along $K$ and  intersecting $\pa\Om$ perpendicularly.
\begin{center}

\bigskip\bigskip


\bigskip

\centerline{\bf AMS subject classification: 53A10, 53C21, 35R35}

\end{center}


\section{Introduction}
Let $\Om$ be an open bounded subset of $\R^{m+1}$, $m\ge 2$, with smooth boundary $\pa\Om$. Recall that the
 \textit{partitioning problem} in $\Om$ consists on finding, for a  given  $0<v<meas\,(\Om)$, a
 critical point of the perimeter functional $\mathcal{P}(\,\cdot\,,\,\Om\,)$ in the class of sets in $\Om$
 that enclose a volume $v$. Here  $\mathcal{P}(\,E\,,\,\Om\,)$ denotes the \textit{perimeter} of $E$ relative
 to $\Om$.

It is clear that whenever such a surface exits will meet $\pa\Om$
orthogonally and will have a constant mean curvature, see Section \ref{ss:FV-a}. In the light
of standard results in geometric measure theory, minimizers do exist
for any given volume and may have various  topologies (see the
survey by A.Ros \cite{A-Ros}). Actually, up to now the complete
description of minimizers have been achieved only in some special
cases, one can see for example \cite{Burg-Kuwert}, \cite{Ritore-Rosales},
\cite{Ros-Verg} and \cite{Stern-Zumbrun}. However, the study of
existence, geometric and topological properties of stationary surfaces
(not necessarily minimizers)  is far from being complete. Let us mention that Gr\"uter-Jost \cite{Gruter-Jost}, have proved the existence of minimal
 discs into convex bodies; while Jost in \cite{Jost} proved the existence of embedded minimal surfaces of higher genus.
   In the particular case of  the free boundary Plateau  problem,  some rather global existence results were obtained
   by M. Struwe in \cite{S-N-pl}, \cite{S-cmc} and \cite{f-b-mi-surf}. In \cite{fall}, the first author proved the
   existence  of surfaces similar to  half spheres surrounding a small volume  near non-degenerate critical points of the mean    curvature of  $\pa\Om$.
Here we are interested in the existence of  families of stationary
sets $E_\e$ for the perimeter functional relative to $\Om$ having
small volume $meas\,E_\e$ proportional to $\e$. Our result generalizes to higher
dimensional sets the one obtained by the first author in
\cite{fall}. Before stating it some preliminaries are needed.
 We denote by $\calV$ the interior normal vector field along $\pa\Om$. For a given smooth set $E\subset\Om$ with
 finite perimeter, let $\Sigma:=\pa E\cap\Om$ satisfy $\pa \Sigma\subset\pa \Om$ and denote by $N$ its exterior
  normal vector field. For a smooth vector field  $X$ in $\R^m$,  the flow of diffeomorphism $\{F_t\}_{t\in(0,t^*)}$
   of $X$ in ${\Om}$ induces a variation $\{E_t=F_t(E)\}_t$ of $E$.
   Set $A(t)=\mathcal{P}(E_t,\Om)$; $V(t)=\textrm{meas}(E_t)$ and
 $$\zeta(p)=\frac{d}{dt}F_t(p)\mid_{t=0}.$$
  It is well known that by the first variation of the perimeter and
  volume functional, one has
\begin{equation}\label{eq:var-air}
A'(0)=-\int_{\Sigma}m H_{\Sigma}\,\langle\zeta,N \rangle \,dA+\oint_{\del\Sigma}\langle\zeta,\bar {N}\rangle\, ds;
\end{equation}
\begin{equation}\label{eq:var-vol}
V'(0)=\int_{\Sigma}\langle\zeta,N\rangle \,dA,
\end{equation}
where $H_{\Sigma}$ is the mean curvature of $\Sigma$, $N$ its
exterior normal vector field and $\bar {N}$ the exterior normal to
$\pa\Sigma$ in $\Sigma$. A variation is called
\textit{normal} if $\zeta=\omega\, N$ for a smooth function
$\omega$, \textit{admissible} if both $ F_t(int \Sigma)\subset \Om$
and $F_t(\pa \Sigma)\subset \pa\Om$ and \textit{volume-preserving}
 if $V(t)=V(0)$ for every $t$. Since for any smooth $\omega$ satisfying $\int_{\Sigma}\omega\,dA=0$
 there exits a volume-preserving admissible normal variation of $E$ with $\zeta=\omega\, N$,
 then
   $E$ is stationary for the perimeter functional ($A'(0)=0$) for any volume-preserving admissible normal
    variation of $E$, if and only if
$$
mH_{\Sigma}\equiv const.\quad\textrm{in } {\Sigma}\qquad\textrm{and}\qquad  \bar{N}(\s)\perp T_{\s}\del\Om\quad\textrm{for every }\s\in
 \pa{\Sigma}.
$$
Up to a change of variable, we can reformulate our question  to the following
free boundary problem: for a given real number $H$, find a
hypersurface $\Sigma\subset\Om_\e$ satisfying the following
conditions
\begin{equation}\label{eq:eq-to-solv}
\begin{array}{ccc}
H_{\Sigma}&\equiv& H\quad\textrm{in } {\Sigma},\\[3mm]
\pa \Sigma&\subset&\pa\Om_\e,\\[3mm]
\la N,\calV^{\e}\ra&=&0\quad\textrm{on } \pa{\Sigma},
\end{array}
\end{equation}
where $\Om_\e:=\e^{-1}\,\Om$ and $\calV^{\e}$ the interior normal vector field on $\pa\Om_\e$.\\
If $K$ is a $k$-dimensional smooth submanifold of $\pa\Om$, we let
$n:=m-k$ and define $K_\e:=\e^{-1}\,K$. Consider the
``half"-geodesic tube contained in $\Om_\e$ around $K_\e$ of radius
1

\[
\bar S_{\e}(K_\e) := \{ q \in \bar{\Om}_{\e}: \quad d (q,K_{\e})=1
\},
\]
with
%
$$
d(q,K_{\e}):=\sqrt{|\mbox{dist}^{\pa\Om_\e}(\tilde{q},K_\e)|^2+|q-\tilde{q}|^2}
$$
where $\tilde{q}$ is  the projection of $q$ on $\pa\Om_\e$ and
$$
\mbox{dist}^{\pa\Om_\e}(\tilde{q},K_\e)=\inf\left\{\mbox{length}(\gamma)\quad:
\quad \gamma\in C^{1}([0,1])\textrm{ is a geodesic in } \pa\Om_{\e}  ;
 \,\,\gamma(0)\in K_\e;\,\,\gamma(1)=\tilde{q}\right\}.
$$
 By the smoothness of $\pa\Om$ and $K$, the tube is a smooth, possibly immersed,
 hypersurface provided $\e$ is sufficiently small.
This tube by construction meets $\pa\Om_\e$ perpendicularly.
Furthermore the mean curvature of this tube satisfies (see also \S~\ref{sss:mc})
\begin{equation}
mH ( \bar S_\e (K_\e) ) = {n}  + \calO (\e)
\label{eq:1-1}
\end{equation}
as $\e$ tends to zero and hence it is plausible under some rather mild assumptions on $K$ that we might be
able to perturb this tube to satisfy (\ref{eq:eq-to-solv}) with  $mH \equiv {n} $. It turns
out that this is not known to  be possible for every (small) $\e > 0$
 but we prove the following theorem~:
\begin{theorem}
\label{th:existence}
Let $\Om$ be a smooth bounded domain of $\R^{m+1}$, $m\ge 2$.
Suppose that $K$ is a non-degenerate minimal submanifold of $\pa\Om$. Then,
there exist a sequence of intervals $I_i = ( \e_i^{-},
 \e_i^{+})$, with $\underline \e_{-} < \overline \e_{+}$
and $\lim_{i \rightarrow +\infty} \e_i^{+} =0$ such that, for all
$\e\in I : =  \cup_i I_i$ the ``half" geodesic tube $\,\e\,\bar S_{\e}(K_\e)\,$
may be perturbed to a hypersurface $\e S_\e$ satisfying (\ref{eq:eq-to-solv})
with mean curvature $ H_{\e S_\e}\equiv \frac{m}{n}\, \e^{-1} $. Namely there
exists a family of  embedded constant mean curvature hypersurface in $\Om$
with boundary on $\pa\Om$ and  intersecting it perpendicularly.
\end{theorem}
\begin{remark}
\begin{itemize}
\item We emphasize that our argument provides also a stationary
area separating of $\R^{m+1}\setminus\bar{\Om}$ when considering
the lower hemisphere parameterized by the
 stereographic projection from the north pole over the unit ball see Section \ref{ss:gt}.
\item  Notice that the surfaces  we obtained might have interesting  topology.
In fact as far as $\e$ tends to zero, our solutions concentrate along $K$ hence  inherit its
 topological structure. Furthermore we cite that some existence result of various  minimal
 immersions were obtained in \cite{Lawson} and   \cite{Schoen-Yau}. \\
We believe that the minimality condition on $K$ should also be
necessary to obtain a result in spirit of Theorem
\ref{th:existence} see the last paragraph of \cite{mp}). The
non-degeneracy condition might fail in some interesting
situations, for example when a
symmetry is present. In this case
however, one can take advantage of it working in a subclass of
invariant functions: this might also guarantee existence for all
small $\e$, see \cite{mp} Section 5.
    \item The hypersurface $S_\e$ is a small perturbation of $\bar
S_\e(K_\e)$ in the sense that it is the normal graph (for some
function whose $L^\infty$ norm is bounded by a constant times $\e$)
over a small translate of $K_\e$ in $\pa\Om_\e$ (by some translation
whose $L^\infty$ norm is bounded by a constant), we refer to Section
\ref{s:appsol} for the precise formulation of the construction of
$S_\e$.
\item This result also remains true for the existence of  \textit{capillary} hypersurfaces in $\Om$
 namely those  with stationary area which intersect $\pa\Om$ in  a constant angle $\gamma\in(0,\pi)$
 along there boundaries. For more precise comments see Remark \ref{rem:ex-cap}.
\end{itemize}
\end{remark}
To prove the theorem, following \cite{mmp}, \cite{mp} and \cite{Ye-1}, we parametrize all
surfaces nearby $\bar S_\e (K_\e)$ having boundaries in $\pa\Om_{\e}$ by two parametric
functions    $\Phi: K\to \R^n$ and $w: S^n_+\times\e^{-1} K\to \R$. Here
$$
S^n_+:=\left\{x=(x^1,\cdots,x^{n+1})\in \R^{n+1}\quad:\quad |x|=1\textrm{ and }x^{n+1}>0\right\}.
$$
This yields a  perturbed tube $S_\e (w,\Phi)$. A standard
computations show that the mean curvature $H(w,\Phi)$ of $S_\e
(w,\Phi)$ is constant, with the right boundary conditions, is
equivalent to solve a system of nonlinear partial differential
equations where the principal part is the Jacobi operator about  a
hypersurface close to $\bar S_\e (K_\e)$. The solvability is based
on the invertibility of this linear operator depending on $\e$
(small parameter). As we will see later, it turns out that this is
possible only for some values of $\e$ tending to zero. Once we
have the invertibility our problem becomes readily a fixed point
problem  that we can solve provided our approximate solution is
accurate enough.
 Our method here is similar in spirit to the
one in \cite{mmp}. It goes back to Malchiodi-Montenegro in
\cite{malm} (see also \cite{mm}, \cite{mal} and \cite{malm2}, for
related issues).\\
 To begin the procedure, we construct first an approximate solution
in the following way: let
   $(y^1,y^{2}\,\dots\,,y^k)\in\R^k$ (resp. $(z^1,z^{2}\,\dots\,,z^n)\in B^n_1$) be the local
   coordinate variables on $K_\e$ (resp. on $S^n_{+}$). Letting
$\Phi:K\to \R^n$ and $w:B^n_1\times K_\e\to \R$,  consider
$$
S_0(y,z)\quad:(y,z)\longrightarrow y\times \e^{-1}\Phi(\e y)\,+ \,(1+w(y,z))\,\T(z).
$$
The nearby surfaces  of $\bar S_\e (K_\e)$ are parameterized
(locally) by
$$
G(y,z):(y,z)\longrightarrow S_0(y,z)\longrightarrow F^\e(S_0(y,z))
$$
where $F^\e$ is  defined in (\ref{eq:par-tu-O}) is ``an almost
isometry" parameterizing a neighborhood of $K_\e$ in $\Om_\e$,
$B^n_1$ is the unit ball centered at the origin and $\T$ is the
stereographic projection from the south pole.
 Call the image of this map
$S_\e(w,\Phi)$, so in particular
\[
S_{\e}(0,0) = \bar S_\e(K_\e).
\]
Notice that since $\Tm{\Big|}_{{_{\pa B^n_1}}}=0$, it follows that
all these surfaces close to $S_\e(K_\e)$ parameterized in this way
have boundaries on $\pa\Om_\e$.

\

Using standard arguments, we compute the mean curvature of
$S_\e(w,\Phi)$, in \S~\ref{sss:mc}.
  The linearized mean curvature operator  about $\bar S_\e (K_\e)$ splits into some
  linear operators on $w$ and $\Phi$, given by
$$
-{\cal L}_\e \, w-  \e \,  \la {\mathfrak J} \, \Phi, \TT\ra +\e
\mathcal{L}^1 w+ \e \mathcal{J}^1(\Phi)+\e^2\,L(w,\Phi),
$$
where ${\mathfrak J}$ is the Jacobi operator about $K$ in the
supporting surface $\pa\Om$, see \S~\ref{ss:Jac-K};
$$
{\cal L}_\e:=\e^2 \, \Delta_K + \Delta_{S^{n}_+} + n;\qquad\mathcal{J}^1\Phi:
=-(3n+1)\,\Tm h(\TT)^a \la \Phi_{\bar{a}},\TT \ra+\,\Tm h(\Phi_{\bar{a}})^a
  +2\Tm h:\Gamma(\Phi)
$$
and $ \mathcal{L}^1$, $L(w,\Phi)$ are  second order differential operators, see \S~\ref{ss:not},
here $ h$ (resp. $\Gamma$) is the second fundamental form of  $\pa\Om$ (resp. $K$) and $h:\Gamma=h_{ab}\Gamma_{ab}$,
where summation over repeated indices is understood.
 The quadratic part of the mean curvature is given by
\begin{eqnarray}\nonumber
&&\frac n2(\e w_{\bar{a}}+\la \Phi_{\bar{a}},\TT \ra)^2
-\e\la\Phi_{\bar{a}},\nabla_{S^n}w_{\bar{a}} \ra-2\e^2 \nabla^2_{K}w:\Gamma(\Phi)\\[3mm]
&&\qquad\qquad+\frac{n+2}{6}\la R(\Phi,\TT)\Phi\,,\,\TT \ra
-\frac 13 \la R(\Phi,E_i)\Phi\,,\,E_i \ra +Q(w)+\e\,Q(w,\Phi).
\end{eqnarray}
Finally the boundary condition reads
\begin{equation*}
\la N,\calV^{\e}\ra=(-1+w)\,\frac{\pa w}{\pa \eta}+\bar{{\cal O} }(\e^2)
+\e^2\, \bar{L}(w,\Phi)+\e\, \bar{Q}(w,\Phi)\quad \hbox{ on } \pa (SNK)_+,
\end{equation*}
where $\eta=-E_{n+1}$ is the normal vector field of $\pa S^n_+$ in
$S^n_+$. \\
The method adopted requires to find an
approximate solution as accurate as possible. For that, letting
$r\geq1$ be an integer and  setting
$$
\hat{w}^{(r)}=\sum_d^r\e^d w^{(d)} \quad\textrm{and} \quad\hat{\Phi}^{(r)}=\sum_d^{r-1} \e^d \Phi^{(d)},
$$
we  solve
\begin{equation}\label{eq:n-perfom}
\begin{array}{ccc}
m\,H(\hat{w}^{(r)},\hat{\Phi}^{(r)})&=& n+\calO(\e^{r+1})\quad\textrm{ on }\quad S_\e
(\hat{w}^{(r)},\hat{\Phi}^{(r)}),\\[3mm]\qquad
 \la N,\mathcal{V}^{\e}\ra&=&\bar{\calO}(\e^{r+2})\quad\textrm{ along }\quad\pa S_\e(\hat{w}^{(r)},\hat{\Phi}^{(r)}).
 \end{array}
\end{equation}
 This leads to an iterative scheme.
 The term of order $\calO(\e)$ appearing in the expansion of the mean curvature (\S~\ref{sss:mc})
 depend linearly on the tangential curvature of  $K$ which is in the kernel of $\Delta_{S^{n}_+} + n$
 and normal curvature $K$ which is perpendicular to this operator. Consequently  by Fredholm theorem,
  we can kill these terms by $w^{(1)}$ provided $K$ is minimal.

  Now to annihilate the higher order
   terms with suitable couples $(w^{(d)},\Phi^{(d-1)})$, $d\geq2$, if we project on the kernel $ \Delta_{S^{n}_+} + n$,
    there
   appears only $\mathfrak{J}$ (the Jacobi operator about $K$) acting on $\Phi^{(d-1)}$ because when we project,
    the term $\mathcal{J}^1\Phi^{(d-1)}$ disappear by oddness.  Moreover neither the nonlinear  terms
     appearing in the expansion  of $H(w,\Phi)$ nor the perpendicularity condition will  influence the iteration as well.
     Therefore  nondegenerency of $K$  is sufficient for this procedure at each step of the iterative scheme.
      In this way for any integer $r\geq1$ we will be able to have (\ref{eq:n-perfom})
      yielding  good approximate solutions. We notice that it is more
      convenient to use the operator $\Delta_{S^{n}_+} + n+\la\mathfrak{J},\TT\ra$ to accomplish this task because
       it is invertible in $L^2(S^n_+\times K)$. Unfortunately one cannot use  it for full solvability of the
        problem because $w$ may not gain regularity. We refer to Section \ref{s:appsol} for more details.
The final step (see \S~\ref{s:span}) is more delicate and consists of the invertibility of the Jacobi  operator about
$S_\e(\hat{w}^{(r)},\hat{\Phi}^{(r)})$ which we call $\mathbb{L}_{\e,r}$.  Let us mention that at this level all
 terms in the expansion depend on $r$ except the model operator $ -{\cal L}_\e \, w-  \e \,
   \la {\mathfrak J} \, \Phi, \TT\ra$. At first glance one sees that  the operator
   $\mathbb{L}_{\e,r}$ is not so close to the model one in the usual Sobolev norms because of the
   competition between the operators $\la {\mathfrak J} \, \Phi, \TT\ra$ and $\mathcal{L}^1_r$.
   This  is due to fact that if one consider  a tube of radius  $\e$ in a manifold $\calM$ with boundary
    sitting on the boundary $\del\calM$, the mean curvature expansion  makes appear terms of order $\e$ depending
    on the second fundamental form of $\del\calM$.  On the contrary,   dealing with manifolds without boundary,
    as in \cite{mmp}, it turns out that in this case the first error terms are of order $\e^2$ and thus also
    in the expansion of the mean curvature of there perturbed tube, there cannot appear terms like $\e L$,
    see \cite{mmp} Proposition 4.1. Having bigger error terms than those in \cite{mmp}, we need  more accurate
    approximate solutions and different spaces  the spectral analysis.
   Since
   our operator  $\mathbb{L}_{\e,r}$ acts on the couple $(w,\Phi)$ almost separately, to tackle this it
   is natural to adjust
    the norms used for $w$ and $\Phi$. For any $v\in L^{2}(S^n_+\times K)$ we decompose it as
    $v=\e^{1-2s}\,w+\la \Phi, \TT\ra$ where $\Phi^i$, $i=1,\dots, n$ are the components of the
     projection of $v$ in the Kernel of $\Delta_{S^{n}_+} + n$ for some $s\in(0,1/2)$. With this
      decomposition, in a suitable weighted Hilbert subspace of $L^{2}(S^n_+\times K)$ we can see
       $\mathbb{L}_{\e,r}$ as a perturbation of the model one, see Proposition \ref{p:model}.\\

 As mentioned above the existence of families of CMC
surfaces only for a suitable sequence  of intervals with length
decreasing to zero and not the whole $\e$ is related to a resonance
phenomenon peculiar to concentration on positive dimensional sets
and it appears in the study of several class of (geometric)
non-linear PDE's. Concentration along sets of dimension $k = 1,
\dots, n-1$ has been proved here, and analogous spectral properties hold
true. By the Weyl's asymptotic formula, if solutions concentrate
along a set of dimension $d$ the average distance between those
close to zero is of order $\e^d$. The resonance phenomenon was taken
care of using a theorem by T. Kato, see \cite{Kato}, page 445, which
allows to differentiate eigenvalues with respect to $\e$. In the
aforementioned papers it was shown that, when varying the parameter
$\e$, the spectral gaps near zero almost do not shrink, and
invertibility can be obtained for a large family of epsilon's. The
case of one dimensional limit sets can be handled using a more direct
method based on a Lyapunov-Schmidt reduction, indeed in this case
the distance between two consecutive small eigenvalues, candidates to
be resonant, is sufficiently large and working away from resonant
modes one can perform a contraction mapping argument quite easily.
Here instead the average distance between two consecutive
eigenvalues becomes denser and denser, to overcome this problem one
needs to apply Kato's Theorem constructing first good approximate
eigenfunctions.

\section{Preliminaries}\label{s:pls}
%
%

Let $K$ be a $k$-dimensional submanifold of $(\partial\Om,\ov g)$
($1\le k\le m-1$) and set $n=m-k$. We choose
along $K$ a local orthonormal frame field $((E_a)_{a=1,\cdots
k},(E_i)_{i=1,\cdots, n})$ which is oriented and call $\mathcal{V}$ the interior normal field
along $\pa\Om$  and $\mathcal{V}_{|K}=E_{n+1}$. At points of $K$, $\R^{m+1}$ splits naturally as
 $T\pa \Om\oplus\R E_{n+1}$ with  $T\pa \Om=T K \oplus N K$, where $T K$ is the
tangent space to $K$ and $N K:=NK^{\pa\Om}$ represents the normal bundle in $\pa\Om$, which
are spanned respectively by $(E_a)_a$ and $(E_j)_j$.
\subsection{Fermi coordinates on $\pa\Om$ near $K$}\label{ss:fc}
Denote by $\n$ the connection induced by the metric $\ov{g}$ and by
$\n^\perp$ the corresponding normal connection on the normal bundle.
Given $q \in K$, we use some geodesic coordinates $\oy$ centered at
$q$.
\begin{equation}\label{eq:fc}
f:\;     {\ov y}\longrightarrow \exp_q^{K}({\ov y}^a E_a).
\end{equation}
 This yields the coordinate vector fields $\ov{X}_a:=f_{*}(\del_{\bar{y}^a})$.
We also assume that at $q$ the normal vectors $(E_i)_i$, $i =
1, \dots, n$, are transported parallely (with respect to $\n^\perp$)
through geodesics from $q$, so in particular
\begin{equation}\label{eq:parall}
    \ov g\left(\nabla_{E_a}E_j\,,E_i\right)=0  \quad \hbox{ at } q,
    \qquad \quad i,j = 1, \dots, n, a = 1, \dots, k.
\end{equation}
In a neighborhood of $q$, we choose {\em Fermi coordinates} $(\oy,
\z)$ on $\pa \Om$ defined by~
\begin{equation}\label{eq:fc}
{\ov F}:\;    ( y,\z)\longrightarrow \exp^{\pa \Om}_{f ({\ov y})}(
\sum\limits_{i=1}^{n}\,\z^i\,E_i); \qquad \quad (\oy, \z) =
\left((\ov y^a)_a,(\z^i)_i\right).
\end{equation}
Hence  we have the coordinate vector fields
$$
\ov{X}_i:=\bar{F}_{*}(\del_{\z^i})\qquad{and}\qquad \ov{X}_a:=\bar{F}_{*}(\del_{\bar{y}^a}).
$$
By our choice of coordinates, on $K$ the metric $\ov{g}_{\al,\be}:=\la \ov{X}_{\al},\ov{X}_{\be}\ra$ splits in
the following way
\begin{equation}\label{eq:splitovg}
    \ov g(q) = \ov g_{ab}(q)\,d\ov{y}^a\otimes d\ov{y}^b+\ov
g_{ij}(q)\,d\z^i\otimes d\z^j; \qquad \quad q \in K.
\end{equation}
We denote by $\Gamma_a^b(\cdot)$ the 1-forms defined on the normal
bundle of $K$ by~
\begin{equation}\label{eq:Gab}
    \Gamma_a^b(E_i)=\ov g(\nabla_{E_a}E_b,E_i).
\end{equation}
We will also denote by $R_{\al\be\gamma\delta}$ the components of the
curvature tensor with lowered indices, which are obtained by means
of the usual ones  $R_{\beta\gamma\delta}^\sigma$ by~
\[R_{\al\be\gamma\delta}=\ov g_{\al\sigma}\,R_{\be\gamma\delta}^\sigma.\]
When we consider the metric coefficients in a neighborhood of $K$,
we obtain a deviation from formula \eqref{eq:splitovg}, which is
expressed by the next lemma, see Proposition 2.1 in \cite{mmp} for
the proof. Denote by $r$ the distance function from $K$.

\begin{lemma}  \label{l:oovg} In the above coordinates $(\oy, \z)$,
for any $a=1,...,k$ and any $i,j=1,...,n$, we
have
\[
\begin{array}{rllll}\ov g_{ij}(0,\z)&=\delta_{ij}+\frac{1}{3}\,R_{istj}\,\z^s\,\z^t\,
+\,{\mathcal O}(r^3);\\[3mm]
\ov g_{aj}(0,\z)&={\mathcal O}(r^2);\\[3mm]
\ov
g_{ab}(0,\z)&=\delta_{ab}-2\,\Gamma_{a}^b(E_i)\,\z^i+\left[R_{sabl}+\Gamma_{a}^c(E_s)\,
\Gamma_{c}^b(E_l) \right]\z^s\z^l+{\mathcal O}(r^3).
\end{array}
\]
Here $R_{istj}$ are  computed at the point $q$ of $K$
parameterized by $(0,0)$.
\end{lemma}
The boundary of the  scaled domain $\de\Om_\e:=\frac{1}{\e}\pa\Om$
is parameterized, in a neighborhood
 of $\e^{-1} q\in K_{\e}:=\e^{-1} K$ by
$$
\bar{F}^{\e}({ y},x'):= \frac{1}{\e}\bar{F}(\e{ y},\e x')\qquad \textrm{ with }x':=(x^i,\cdots,x^n).
$$
Hence  we have the induced coordinate vector fields
$$
{X}_i:=\bar{F}^{\e}_{*}(\del_{x^i})\qquad\textrm{and}\qquad {X}_a:=\bar{F}^{\e}_{*}(\del_{{y}^a}).
$$
 By construction, ${X_{\al}}_{|\e^{-1} q}=E_{\al}$ and  $\mathcal{V}^{\e}(\e^{-1} q)=E_{n+1}$.
From Lemma \ref{l:oovg} it is evident that the metric $g$ on $(\pa\Om_\e,g)$ has the expansion given by the
\begin{lemma}\label{l:ovg}
In a neighborhood of $K_\e$ the following hold
\[
\begin{array}{rllll} g_{ij}(0,x)&=\delta_{ij}+\frac{\e}{3}\,R_{istj}\,x^s\,x^t\,
+\,{\mathcal O}(\e^2 r^3);\\[3mm]
 g_{aj}(0,x)&={\mathcal O}(\e r^2);\\[3mm]

g_{ab}(0,x)&=\delta_{ab}-2\,\Gamma_{a}^b(E_i)\,x^i+\e\,\left[R_{sabl}+\Gamma_{a}^c(E_s)\,
\Gamma_{c}^b(E_l) \right]x^sx^l+{\mathcal O}(\e^2r^3).
\end{array}
\]
\end{lemma}
We can now    parameterize  tubular neighborhood of $K_\e$ in
$\Om_{\e}$,
\begin{equation}\label{eq:par-tu-O}
F^{\e}({ y},x',x^{n+1})=\frac{1}{\e}\bar{F}(\e{ y},\e x')+x^{n+1}\mathcal{V}^{\e}({ y},x'),
\end{equation}
 where $\mathcal{V}^{\e}({ y}, x'):=\mathcal{V}(\frac1\e\bar{F}(\e{ y},\e x'))$.
 We denote  by $h$ the second fundamental form of $\pa\Om$  so that:
\begin{equation}\label{eq:scn-ffm}
\langle d \mathcal{V}^{\e}(p)[X_{\alpha}],X_{\beta}\rangle=\e\,h_{\alpha,\beta}(q)
\end{equation}
when $q={\bar{F}^{\e}}(p)$.
\subsection{The Jacobi operator about $K$}\label{ss:Jac-K}
The linearized mean curvature operator about $K$ is given by
\begin{equation}
{\mathfrak J} : = \Delta^{\perp} - {\cal R}^{\perp} + {\mathcal B}
\end{equation}
where the normal Laplacian $\Delta^{\perp}$ is defined as \[
\Delta^{\perp} := \nabla^{\perp}_{E_a} \, \nabla^{\perp}_{E_a} -
\nabla^\perp_{\nabla^T_{E_a}E_a},
\]
with $\n^\perp$ denoting the connection on the normal bundle of $K$ in $\pa \Om$. While
${\mathcal B}$ is a symmetric operator
defined by
\[
\la {\mathcal B} (X), Y \ra  = \Gamma_a^b (X)\, \Gamma_b^a (Y)\qquad \hbox{for all } X,Y\in T_pK,
\]
where $\Gamma$ is defined in \eqref{eq:Gab} and ${\cal R}^\perp : N_p K\longrightarrow N_p
K$ is given  by
\[
{\cal R}^\perp  : = \left( R(E_a, \cdot )\, E_a \right)^\perp,
\]
where $( \cdot )^{\perp}$ denotes the orthogonal projection on $N_p K$.
Finally, we recall that the Ricci tensor is defined by
\[
\mbox{Ric}  (X, Y) = -   \la R (X, E_\gamma)  \, Y, E_\gamma\ra\qquad \hbox{for all } X, Y \in T_p M.
\]
\subsection{First and second variation of area for capillary hypersurfaces}\label{s:FSV-a}
Let $\Sigma$ be a smooth hypersurface in an $(m+1)$-dimensional Riemannian manifold $(M,g)$ with smooth,
nonempty boundary $\pa M$.
Suppose that $\del\Sigma\subset\del M$ so that $M$ is separated into two parts, call $\Lambda$ the boundary of
one of these parts in $\pa M$.
\subsubsection{First variation of area}\label{ss:FV-a}
  Let $F_t$  be a variation of $\Sigma$ with variation  vector field
  $$
  \z(p)=\frac{\del F_t}{\del t}(p)_{|t=0}\quad\textrm{ for every }p\in \Sigma.
  $$
A variation is called \textit{admissible} if both $ F_t(int \Sigma)\subset M$ and $F_t(\del \Sigma)\subset \del M$.
Let $N$ be a unit normal vector along $\Sigma$; $H_{\Sigma}$ its mean curvature
 and $\upsilon$ (respectively $\bar{\upsilon}$) be the unit exterior normal vector along $\del \Sigma$ in $\Sigma$
 (respectively in $\Lambda$).

 An admissible variation  induces  hypersurfaces $\Sigma_t$ and $\Lambda_t$. Let  $A(t)$ (respectively $T(t)$) be
 be the volume of $\Sigma_t$ (respectively $\Lambda_t$) and $V(t)$ the signed volume bounded by $\Sigma$ and $\Sigma_t$.
 For a given angle $\gamma\in(0,\pi)$, we consider the total energy
 \begin{equation}  \label{eq:tot-nrj}
 \calE (t):=A(t)-\cos(\gamma)\, T(t).
 \end{equation}
It is well known  (see for example  \cite{Ros-Souam}) that
\begin{equation}\label{eq:var-air}
\calE'(0)=-\int_{\Sigma}n H_{\Sigma}\langle\zeta,N \rangle_g dA+\oint_{\del\Sigma}\langle\zeta,\upsilon-
\cos(\gamma)\,\bar{\upsilon}\rangle_g ds
\end{equation}
and
\begin{equation}\label{eq:var-vol}
V'(0)=\int_{\Sigma}\langle\zeta,N\rangle_g dA.
\end{equation}
A variation is called \textit{volume-preserving} if $V(t)=V(0)$ for every $t$.
$\Sigma$ is called \textit{ capillary hypersurface} if $\Sigma$ is stationary for
 the total energy ($\calE'(0)=0$) for any volume-preserving admissible  variation.
 Consequently  if $\Sigma$ is capillary, it has a constant mean curvature and intersect
 $\pa M$ with the angle $\gamma$ in the sense that  the angle between the normals of $\upsilon$
  and $\bar{\upsilon}$ is $\gamma$ or equivalently  the angle between   $N$ and $\calV$ is $\gamma$,
   where $\calV$ is the unit outer normal field along $\pa M$.\\
Physically, in the tree-phase system the quantity $\cos(\gamma)\,T(0)$  is interpreted as the
\textit{wetting energy} and $\gamma$ the \textit{contact angle} while  $\cos(\gamma)$ is the
\textit{relative adhesion coefficient} between the fluid bounded by $\Sigma$ and $\Gamma$ and
 the walls $\pa M$. Here we are interested in a configuration in the absence of gravity. A more
  general setting  including the gravitational energy and works on capillary surfaces can be found
  in the book by R. Finn \cite{Finn}.

%
\subsubsection{ The Jacobi operator about $\Sigma$}\label{ss:SV-a}
We denote by $\Pi_{\Sigma}$ and $\Pi_{\del M}$ the second  fundamental form of $\Sigma$ and of
 $\del M$ respectively. Assume that $\Sigma$ is a capillary hypersurface. Recall that the Jacobi
  operator (the linearized mean curvature operator about $\Sigma$) is given by the second variation
  of the total energy functional $\calE$. For any volume-preserving admissible  variation, we have
  (see \cite{Ros-Souam} Appendix for the proof)
\begin{equation}\label{eq:secvar-air}
\calE''(0)=-\int_{\Sigma}\left(\omega\Delta_{\Sigma} \omega+|\Pi_{\Sigma}|^2 \omega^2+Ric_g(N,N)\omega^2\right)
dA+\oint_{\del\Sigma}\big(\omega\frac{\del \omega}{\del\upsilon}-q\,\omega^2\big)ds,
\end{equation}
 where
$$
\omega=\la\zeta,\,N\ra_g \qquad \textrm{and}\qquad q=\frac{1}{\sin(\gamma)}\,\Pi_{\pa M}
(\bar{\upsilon},\bar{\upsilon})-\cot(\gamma)\,\Pi_{\Sigma}({\upsilon},{\upsilon}).
$$
Since for any smooth $\omega$ with $\int_{\Sigma}\omega dA=0$ there exits an admissible,
volume-preserving variation with variation vector field $\omega\, N$ as a normal part,
we have now the Jacobi operator about $\Sigma$ that we define by duality as
$$
\la \mathfrak{L}_{\Sigma,N}\,\omega,\omega'\ra:=\int_{\Sigma}\left\{\n\omega\n \omega'
-\left(|\Pi_{\Sigma}|^2+Ric_g(N,N)\right)\omega\,\omega'\right\}dA+\oint_{\del\Sigma}q\,\omega\,\omega'ds.
$$
%
%
%
\begin{remark}\label{rem:linkJac}
Let us  observe that for any smooth $\hat{\omega}\,N$ and  $\hat{N }$  transverse vector
field  along $\Sigma$ there induce an admissible volume preserving variation. The linearized
 mean curvature operators $\mathfrak{L}_{\Sigma,N}$ and $\mathfrak{L}_{\Sigma,\hat{N}}$ are linked by
$$
\mathfrak{L}_{\Sigma,\hat{N}}\,\hat{\omega}=\mathfrak{L}_{\Sigma,N}\,({\la N,\hat{N}\ra_g\,\hat{\omega}})
+m\,\hat{N}^T(H_{\Sigma})\,\hat{\omega},
$$
where $\hat{N}^T$ is the orthogonal projection of $\hat{N}$ on $T\Sigma$. This shows that
$\mathfrak{L}_{\Sigma,\hat{N}}$ is self-adjoint with respect to the inner product
$$
\int_{\Sigma}\hat{\omega}\,\hat{\omega}'\,\la N,\hat{N}\ra_g\,dA.
$$
\end{remark}
\subsection{The stereographic projection}\label{ss:stgp}
We will denote by $\bfp:\R^n\to S^n$  the inverse of the stereographic projection from the south pole.
 $\bfp=\left(\,\bfp^1\,,\dots,\,\bfp^n,\,\bfp^{n+1}\,  \right)$ is a conformal parametrization of
 $S^n$ and for any $z=(z^1,\dots ,z^n)\in \R^n$,
\begin{eqnarray*}
\bfp(z)&=&(z,1)\,\mu(z)-E_{n+1}\\[3mm]
     &=&\left(\,\frac{2\,z^1}{1+|z|^2}\,,\dots,\,\frac{2\,z^n}{1+|z|^2},\,\frac{1-|z|^2}{1+|z|^2}\right)\\[3mm]
\end{eqnarray*}
with conformal factor given by
\begin{equation}\label{eq:mu}
\mu(z):=\frac{2}{1+|z|^2}.
\end{equation}

 We often use the projection of $\bfp$ on $\R^n$ and denote it by
\begin{equation}\label{eq:tibfp}
\tilde{\bfp}(z):=(z,0)\,\mu(z).
\end{equation}
We collect in the following lemma some properties of the function $\bfp$
 which will be useful later on, we omit here the proof which can be obtained rather easily with elementary
 computations
\begin{lemma}\label{l:formulas} For every $i,j,l=1,\dots,n$, there holds
$$\langle \bfp_{i},\bfp_j \rangle=\mu^2\,\de_{ij}; \qquad\bfp_i^{n+1}=-\mu \,\bfp^i;\qquad
\tilde{\bfp}_i=-\bfp^i\,\tilde{\bfp}+\mu\, E_i;$$
$$
\langle \bfp_{ii},\bfp_l \rangle=\mu^2\,\bfp^l-2\mu^2\,\bfp^i\, \de_{il}.
$$
\end{lemma}
Recall that the Laplace operator on $S^n$ can be expressed in terms of the Euclidean one by the formula
$$
\Delta_{S^n}=\frac{1}{\mu^2}\left( \Delta_{\R^n}+(2-n)\bfp^i \pa_{i}  \right).
$$
Moreover, it is easy to verify that
$$
\Delta_{S^n}\bfp+n\bfp=0.
$$
 It is clear that for any $0<r\leq1$ the restriction of $\bfp$ on $B^n_r$ parametrizes
  a spherical cap $S^n(r)$, where $B^n_r$ is  a the ball centered at $0$ with radius $r$. \\
 Given $\gamma\in(0,\pi)$, if we let $r^2=\frac{1-\cos(\gamma)}{1+\cos(\gamma)}$, the
 image by $\bfp$ of $B^n_r$ is the spherical cap $S^n(\gamma)$ which intersects the horizontal
  plane $\R^n+\cos(\gamma)\,E_{n+1}$ and   makes an angle $\gamma$ with it. In particular we
   denote (henceforth define)
$$
\T(\gamma):=\bfp{\Big|_{B^n_{r(\gamma)}}}-\cos(\gamma)\,E_{n+1};\qquad \T:=\T({\frac{\pi}{2}})
$$
$$
S^n_+:=S^n(\frac{\pi}{2})=\left\{x=(x^1,\dots,x^{n+1})\in \R^{n+1}\quad:\quad |x|=1\textrm{ and }x^{n+1}>0\right\}.
$$

For any $0<r\leq 1$, denote by $\tau_r$ the unit outer normal vector of $\pa B^n_r$,
the normal field (not unitary) of $\pa S^n(r)$ in $S^n(r)$ expressed as follows
$$
\frac{\pa\bfp}{\pa\tau_r}{\Bigg|_{\pa B^n_r} }=\mu\,|\tilde{\bfp}|\,\left( \bfp^{n+1}
\frac{\tilde{\bfp}}{|\tilde{\bfp}|^2}-E_{n+1}  \right){\Bigg|_{\pa B^n_r} }.
$$
Now when $r^2=\frac{1-\cos(\gamma)}{1+\cos(\gamma)}$, the unit normal
in $S^n(\gamma)$ of $\pa S^n(\gamma)$ is given and denoted by
\begin{equation}\label{eq:noSn}
\eta(\gamma)=\cot(\gamma)\,\TT({\gamma})-\sin(\gamma)\,E_{n+1},
\qquad\textrm{in particular}\qquad\eta:=\eta({\frac{\pi}{2}})=-E_{n+1}
\end{equation}
while the unit normal of $\pa S^n(\gamma)$ in the plane $\R^n+\cos(\gamma)\,E_{n+1}$ is
 $\frac{\TT({\gamma})}{|\TT({\gamma})|}|_{\pa B^n_r}$.\\
Observe that the angle between the two normals $\frac{\TT({\gamma})}{|\TT({\gamma})|}$
and $\eta(\gamma)$ is $\gamma$ along $\pa S^n(\gamma)$, namely since ${|\TT({\gamma})|}=\sin(\gamma)$ on $\pa B^n_r$,
$$
\la \frac{\TT({\gamma})}{|\TT({\gamma})|},\eta(\gamma)  \ra=cos(\gamma)\qquad\textrm{on }\pa S^n(\gamma) .
$$
Consider  the eigenvalue problem, $u:S^n(\gamma)\to\R$,
\begin{equation*}
\left\{\begin{array}{ccc}
\displaystyle\Delta_{S^n(\gamma)}u+n u &=& 0 \quad\textrm{  in }S^n(\gamma);\\[3mm]
\displaystyle\frac{\pa u}{\pa\eta(\gamma)}&=&\cot(\gamma)\,u\quad\textrm{ on } \pa S^n(\gamma).
\end{array}\right.
\end{equation*}
 It is well known that  the only solutions to the interior equation  are the degree one
  homogeneous polynomials on $S^n_+$, spanned by the $n+1$  components of $\bfp$. By \eqref{eq:noSn}
   the boundary condition is satisfied only by   $\Theta^i(\gamma)$, $i=1,\cdots, n$.


\subsection{ Notations}\label{ss:not}
 In the following, expressions of the form
$L(w,\Phi)$ denote linear operators, in the functions $w$ and
$\Phi^j$ as well as their derivatives with respect to the vector
fields $\e \, X_{{a}}$ and $X_i$ up to second order, the coefficients
of which are smooth functions on $S^{n}_+ \times K$
bounded by a constant independent of $\e$ in the ${\cal C}^\infty$
topology (where derivatives are taken using the vector fields
$X_{\bar{a}}$ and $X_i$). Also  $\bar{L}(w,\Phi)$ are restrictions of expressions
 like $L(w,\Phi)$ on $\pa S^{n}_+ \times K$ with  $L(w,\Phi)$ contains only one derivative
  of $w$ or $\Phi$ with respect to the vector
fields $\e \, X_{{a}}$ and $X_i$.

\medskip

Similarly, expressions of the form $Q (w,\Phi)$ denote nonlinear
operators, in the functions $w$ and $\Phi^j$ as well as their
derivatives with respect to the vector fields $\e \, X_{{a}}$ and
$X_i$ still up to second order, whose coefficients of the Taylor expansion
are smooth functions on $S^{n}_+ \times K$ which are bounded by a
constant independent of $\e$ in ${\cal C}^\infty$ topology
(where derivatives are taken using the vector fields $X_{{a}}$ and
$X_i$). Moreover, $Q$ vanish quadratically in the pair $(w,
\Phi)$ at $0$ (that is, its Taylor expansion does not involve any
constant nor any linear term). Also  $\bar{Q}(w,\Phi)$ are restrictions
of expressions like $Q(w,\Phi)$ on $\pa S^{n}_+ \times K$ with  $Q(w,\Phi)$
contains only one derivative of $w$ or $\Phi$ with respect to the vector
fields $\e \, X_{{a}}$ and $X_i$.

\medskip

Finally, terms denoted $\calO(\e^d)$ are smooth functions on
$S^{n}_+\times K_\e$ which are bounded by a constant times $\e^d$
in ${\cal C}^\infty$ topology (where derivatives are taken using
the vector fields $X_a$ and $X_i$). Also expressions like $\bar{\calO}(\e^d)$ are restrictions
of $\calO(\e^d)$ on $\pa S^{n}_+ \times K$.

\section{Geometry of tubes}\label{ss:gt}
We derive expansions as $\e$ tends to $0$ for the metric, second
fundamental form and mean curvature of  $\bar S_\e (K_\e)$
and their perturbations.

\subsubsection{Perturbed tubes}

We now describe a suitable class of deformations of the geodesic tubes (in the metric induced by $F^\e$ on $\R^{m+1}$)
 $\bar S_\e(K_\e)$, depending on a section $\Phi$ of $NK_\e$  and
a scalar function $w$ on the spherical normal bundle $(SNK_\e)_+$ in $\pa \Om_\e$.\\
We recall that $(y^1,y^{2}\,\dots\,,y^k)\in\R^k$ (resp. $(z^1,z^{2}\,\dots\,,z^n)\in B^n_1$)
 are the local coordinate variables on $K_\e$ (resp. on $S^n_{+}$). Letting $\Phi:K\to \R^n$
 and $w:B^n_1\times K_\e\to \R$,  consider
$$
S_0\,:(y,z)\mapsto  y\times \e^{-1}\Phi(\e y)\,+ \,(1+w(y,z))\,\T(z).
$$
The nearby surfaces  of $\bar S_\e (K_\e)$ is parametrized (locally) by
$$
G(y,z):(y,z)\longrightarrow S_0(y,z)\longrightarrow F^\e(S_0(y,z))
$$
namely
\[
G(y,z) :=  F^{\e} \, \left( { y},\frac{1}{\e}\Phi(\e y)+(1+w(y,z))\tilde{\T}(z),(1+w(y,z))\T^{n+1}(z)\right).
\]
Since $\Tm{\Big|}_{{_{\pa B^n_1}}}=0$, it follows
$$
G(y,z){\Big|}_{{_{\pa B^n_1}}}\in\pa\Om_\e\qquad\textrm{ for any } y .
$$
%
 The image of this map will be called
$S_\e(w,\Phi)$. In particular
\[
S_{\e}(0,0) = \bar S_\e(K_\e).
\]
It will be understood that for any fixed point $p=F^\e(0,0)\in K_\e$,
$\Phi(\e\,y)\in  NK_\e\subset T_{p}\pa\Om_\e$ and  $\T(z)\in S^n_+\subset NK_\e\oplus \R E_{n+1}$
 are in the tangent space at $p$ of $\R^{m+1}$ endowed with the metric induced by $F^\e$. For more
  convenience we introduce the following  notations

\

\noindent {\bf Notation:} On $K_\e$ we will consider
\[
\Phi  : = \Phi^j \, E_j \qquad \qquad  \Phi_a : = \partial_{y^a}
\, \Phi^j \, E_j \qquad \qquad \Phi_{ab} : = \partial_{y^a}
\partial_{y^b} \, \Phi^j \, E_j
\]
\[
\Theta  : = \Theta^j \, E_j +\Tm E_{n+1}=\TT+\Tm E_{n+1} \qquad  \Theta_i  :=
\partial_{z^i} \Theta^j \, E_j+\partial_{z^i} \Tm E_{n+1}=\TT_i+\partial_{z^i} \Tm E_{n+1}.\\
\]

For simplicity, we will write
\[
w_j :=  \del_{z^j} w ; \quad  w_{a} :=  \del_{ y^a} w ;
\qquad w_{ij} : =  \del_{z^i}\, \del_{z^j} w; \quad  w_{ a
b} : =  \del_{ y^a}\, \del_{ y^b} w ; \quad w_{ a
j} : = \del_{ y^a}\, \del_{z^j} w;\\
\]
It is easy to see that the tangent space to $S_\e(w, \Phi)$ is spanned by the vector
fields
\begin{equation}
\begin{array}{rcccl}
Z_{a} & = & G_* (\del_{ y^a}) & = &  X_a +
w_{ a} \, \Upsilon + \Psi_a +(1+w)\Tm D_a\mathcal{V}^{\e}, \qquad a = 1, \ldots, k
\\[3mm]
Z_j &= &  G_* (\del_{z^j}) & = & (1+w)\, \Upsilon_j +
w_j \, \Upsilon+(1+w)\Tm D_j\mathcal{V}^{\e} , \qquad j=1, \ldots, n,
\end{array}
\label{eq:defz0zj}
\end{equation}
where
\[
\Psi  : = \Phi^j \, X_j; \qquad \qquad  \Psi_a : = \partial_{y^a}
\, \Phi^j \, X_j;
\]
\[
\Upsilon : = \Theta^j \, X_j +\Tm \mathcal{V}^{\e}; \qquad \qquad \Upsilon_i : =
\partial_{z^i} \Theta^j \, X_j+\pa_{z^j}\Tm \mathcal{V}^{\e}
\]

and
\begin{equation}\label{eq:dandaj}
\begin{array}{ccc}
D_a{\mathcal{V}^{\e}( y,(1+w(y,z))
\TT+\e^{-1}\Phi(\e y))}&=&\e\left(h_{a\al}+(w_a \T^l+\Phi_a^l)h_{l\al}\right)X_\al;\\[3mm]

 D_j{\mathcal{V}^{\e}( y,(1+w(y,z))
\TT+\e^{-1}\Phi(\e y))}&=&\e\left(w_j \T^l+(1+w)\T^l_j\right)h_{l\al}X_\al.
 \end{array}
\end{equation}
\subsubsection{The first fundamental form}\label{sss:1ff}
In this subsection we expand the coefficients of the first
fundamental form of $S_\e (w, \Phi)$.  Using the expansions in
 Lemma \ref{l:ovg}, one can easily get
\begin{equation}
\begin{array}{rcl}
\la X_a, X_b\ra  & = & \delta_{ab} - 2 \, \e \, \Gamma^b_{a}
(\Theta) - 2 \, \Gamma^b_{a} \, (\Phi) + {\cal O} (\e^2) +
\e \, L (w, \Phi) + Q (w, \Phi) \\[3mm]
\la  X_i, X_j\ra  & = &  \delta_{ij} + \frac{\e}{3} \big( \la R ( \Theta , E_i) \, \Phi ,
E_j \ra + \la R(\Phi , E_i ) \, \Theta , E_j \ra \big)+ \frac{1}{3} \la R ( \Phi , E_i) \, \Phi ,
E_j \ra\\[3mm]
& + &  \calO (\e^2) + \e^2\,
L (w,\Phi) + \e\,Q (w,\Phi) \\[3mm]
\la X_i, X_a \ra  & = &  {\cal O} (\e^2) + \e \, L(w, \Phi) +
Q(w,\Phi).
\end{array}
\label{lems}
\end{equation}

These together with the fact that $R(\TT , \TT) =0$ imply
\begin{equation}
\la \Upsilon, \Upsilon_j \ra  = \frac{\e}{3} \, \la  R(\Phi,\TT )
\, \TT , \TT_j \ra  + \frac{1}{3} \, \la  R(\Phi,\TT )\, \Phi , \TT_j \ra
+\mathcal{O}(\e^2) +\e^2\, L (w,\Phi) + \e\,Q (w,\Phi)
\label{eq:jjjjj}
\end{equation}

Using similar arguments, and the fact that $\la\Upsilon, \Upsilon\ra  =1$ on $K_\e$ yields
\begin{equation}
\la  \Upsilon, \Upsilon \ra  = 1 +\frac{1}{3} \, \la  R(\Phi,\TT )
\, \Phi , \TT \ra + \mathcal{O}(\e^2) +\e^2\, L (w,\Phi) + Q (w,\Phi) \label{eq:edf}
\end{equation}
Moreover
\begin{eqnarray}
\la \Upsilon_i, \Upsilon_j \ra &=& \la \T_i,\T_j \ra+\frac{1}{3} \left( \la R(\Phi,\TT_i )
\, \TT , \TT_j \ra +\la R(\Phi,\TT_j )
\, \TT , \TT_i \ra \right)\\
&+& \frac{1}{3} \, \la R(\Phi,\TT_i )
\, \Phi , \TT_j \ra +\mathcal{O}(\e^2) +\e^2\, L (w,\Phi) + Q (w,\Phi). \nonumber
\label{eq:UiUj}
\end{eqnarray}
Now, by \eqref{eq:dandaj} we have that
\begin{equation}\label{eq:djnU}
\la D_j{\mathcal{V}^{\e}},\Upsilon \ra = \e (1+w)\la h(\TT),\TT_j \ra+\e w_j \la h(\TT),\TT \ra
+\mathcal{O}(\e^2) +\e^2\, L (w,\Phi) + \e Q (w,\Phi)
\end{equation}
and
\begin{equation}\label{eq:djnUi}
\la D_j{\mathcal{V}^{\e}},\Upsilon_i \ra = \e (1+w)\la h(\TT_i),\TT_j \ra+\e w_j \la h(\TT),\TT_i \ra
+\mathcal{O}(\e^2) +\e^2\, L (w,\Phi) + \e Q (w,\Phi)
\end{equation}
We are now in position to expand the coefficients of the
first fundamental form of $S_\e(w,\Phi)$. We have
\begin{proposition}\label{pr:metric}
For any $a,b \in \{1,\cdots, k\}$ and $i,j\in \{1,\cdots, n\}$, we have that
\begin{equation}\label{eq:zazb}
  \langle Z_{a}, Z_{ b}\rangle  =  \delta_{ab}+2\e\Tm h_{ab}  - 2\e
\Gamma_{a}^b (\TT ) -2 \,
\Gamma_{a}^b (\Phi )+ \calO(\e^2)  + \e \, L(w,\Phi) + Q(w,\Phi)
\end{equation}
\begin{eqnarray}\label{eq:zazj}
  \la Z_{ a} , Z_j \ra & = & 2\e \Tm h(\TT_j)^a+\la\Phi_{\bar{a}},\TT_j \ra
  +\calO(\e^2) + \e L(w,\Phi)+Q(w,\Phi)
\end{eqnarray}
\begin{eqnarray}\label{eq:zizj}\nonumber
 \la Z_i, Z_j \ra & = & \la\Theta_i , \Theta_j \ra\, (1 + 2w) +
   2\e(1+3 w)\Tm \la h(\TT_i),\TT_j \ra \nonumber\\
   & + & 2\e \Tm\left(\la h(\TT_i),\TT \ra w_j
   + \la h(\TT_j),\TT\ra w_i  \right) \nonumber
  \\ & + & \frac{\e}{3}\, \left( \la R(\TT , \TT_i)\, \Phi , \TT_j\ra +
   \la R(\TT , \TT_j)\, \Phi , \TT_i\ra \right)+w_i w_j +\la\Theta_i , \Theta_j \ra w^{(2)}  \\
   & + &
   \frac{1}{3}\, \la R(\Phi , \TT_i)\, \Phi , \TT_j\ra +
\calO(\e^2) +  \e^2 \, L (w,\Phi) + \e Q (w,\Phi).\nonumber
\end{eqnarray}
\end{proposition}

\subsubsection{The normal vector field}\label{sss:nvf}

In this subsection we expand the unit normal to $S_\e
(w,\Phi)$.
Define the vector field
\[
\tilde N  : =  - \, \Upsilon + \al^j \, Z_{j} +  \be^c \,Z_{ c},
\]
it is the outer normal field along $S_\e
(w,\Phi)$ if we can determine  $\al^j$ and $\be^c$ so that $\tilde
N$ is orthogonal to all of the $Z_b$ and $Z_i$. This leads
to a linear system for $\al^j$ and $\be^a$.

\medskip

We have the following expansions
\begin{equation}\label{eq:Uza}
\la  \Upsilon , Z_{ a} \ra  =   w_{ a} +
  \la\Phi_{\bar a} , \TT \ra +\e \Tm \, (h(\TT))^a + \e^2 \,  L(w, \Phi)
  + \e \, Q(w, \Phi);
  \end{equation}
\begin{eqnarray}\label{eq:Uzj}  \nonumber
\la  \Upsilon, Z_j \ra   & = &  w_j +
  \e (1+2 w) \Tm \la h(\TT),\TT_j \ra +2\e \Tm w_j \la h(\TT),\TT \ra
\\
& + & \frac{\e}{3} \, \la
R( \Phi  , \TT ) \, \TT , \TT_j \ra + \frac{1}{3} \, \la
R( \Phi  , \TT ) \, \Phi , \TT_j \ra+\mathcal{O}(\e^2)
+ \e^2 \, L(w,
\Phi) + \e \, Q (w, \Phi),
\end{eqnarray}
 These follow from (\ref{lems})
together with the fact that $ \la\Upsilon , Z_{ a} \ra =0$ and $
\la \Upsilon, Z_j \ra = 0$ on $K_\e$.

\medskip

Using Proposition~\ref{pr:metric}, and some algebraic calculations, one can obtain
\begin{equation}\label{eq:bc}
\be^c = w_{c}  + \la \Phi_{{c}}, \TT \ra + \e \Tm h(\TT)^c+\mathcal{O}(\e^2) +\e \, L(w, \Phi) +
 Q(w, \Phi).
\end{equation}

and
\begin{eqnarray}\label{eq:alj}\nonumber
\al^j \, \la \T_j, \T_i \ra & = &  w_i+ \e \Tm \la h(\TT),\TT_i \ra
   +\e \Tm \la h(\TT),\TT  \ra w_i\\
   & -  & 2\e\Tm \left(\la h(\TT_l),\TT_i  \ra w_l
   + h(\TT_i)^a w_a + h(\TT_i)^a \la \Phi_a,\TT \ra\right)\nonumber\\
   & + & \frac{1}{3} \,\e  \la R (\Phi, \TT) \, \TT , \TT_i \ra
   - \e \Tm h(\TT)^a \la \Phi_a,\TT_i \ra \\
   & - & 2 w w_i   - w_a \la \Phi_a,\TT_i \ra - \la \Phi_a,\TT \ra \la \Phi_a,\TT_i \ra
   + \frac{1}{3} \, \la R (\Phi, \TT) \, \Phi , \TT_i \ra\nonumber \\
   & + & \mathcal{O}(\e^2)+
\e^2 \, L(w, \Phi) + \e Q(w, \Phi).\nonumber
\end{eqnarray}
Using these and the fact that $\la \T_j, \T_i \ra =\mu^2 \de_{ij}$,
a straightforward computations imply
\begin{eqnarray*}
|\tilde N|^{-1} & = & 1 + \e \Tm
   \left( \frac{1}{\mu^2} \la h(\TT),\TT_i \ra w_i +  h(\TT)^c w_c + h(\TT)^c
   \la \Phi_{ c}, \TT \ra \right) +\frac{1}{6} \,  \la R (\Phi, \TT) \, \Phi , \TT \ra\\
  & +& \frac12 \left( w_c^2+\frac{1}{\mu^2} w_j^2  +2 w_c \la \Phi_{c},\TT \ra +
  \la \Phi_{ c},\TT \ra^2  \right)+ \mathcal{O}(\e^2)+
\e^2 \, L(w, \Phi) + \e Q(w, \Phi).
\end{eqnarray*}
The unit normal to the perturbed geodesic tube is then given simply by
$N=\frac{\tilde{N}}{|\tilde{N}|}$. We summarize this in the following lemma
\begin{proposition}\label{pr:normal}
The normal vector field $N$ to $S_\e(w, \Phi)$ is given by $N=\frac{\tilde{N}}{|\tilde{N}|}$
where
\begin{equation}
\tilde N  : =  - \, \Upsilon +  \alpha^j \, Z_j + \beta^c \,Z_c
\label{eq:3-5}
\end{equation}
and where the coefficients $\alpha^j$ and $\be^c$ are given by formulas \eqref{eq:alj} and \eqref{eq:bc}.
\end{proposition}
Using the fact that $\Tm{\Big|}_{{_{\pa B^n_1}}}=0$ we can easily deduce
\begin{lemma}\label{lem:perp}
The perpendicularity condition is given by
$$
\la N,\calV^\e \ra=(-1+w)\, w_j z^j+\bar{{\cal O} }(\e^2)+\e^2\, \bar{L}(w,\Phi)+
\e\, \bar{Q}(w,\Phi)\quad \hbox{ on } \pa (SNK)_+,
$$
%
%
\end{lemma}
{\bf Proof~:}
Since  $\Tm{\Big|}_{{_{\pa B^n_1}}}=0$ it follows that
$\la \calV^\e,-\Upsilon+\be^c\,Z_c \ra=0$ on $\pa B^n_1$ on the other hand using the fact
 that $R(E_i,E_i)=0$ with $\frac{\pa \TT}{\pa\tau}{\Big|}_{{_{\pa B^n_1}}}=0$ (see \S \ref{ss:stgp}) we get
$$
\la\al^j\,Z_j,\calV^\e\ra=(-1+w)\, w_j\T^{n+1}_j+\bar{{\cal O} }(\e^2)+\e^2\, \bar{L}(w,\Phi)+\e\,
 \bar{Q}(w,\Phi)\quad \hbox{ on } \pa (SNK)_+.
$$
The lemma now follows since  $\T_j^{n+1}=-\mu \T^j=-\mu^2 z^j$  and $\mu{\Big|}_{{_{\pa B^n_1}}}=1$.
\hfill $\Box$
%


\subsubsection{The second fundamental form}\label{sss:2ff}
In this subsection we expand the coefficients of the second fundamental form. Recall that
 $\n$ is the Levi-Civita connection on $\pa \Om$ and $h$ its second fundamental form, the
  derivation for vector fields on $\pa\Om$   yields
$$
\frac{\pa}{\pa z^i} X_\al(  y,(1+w(y,z)) \TT+\e^{-1}\Phi(\e
y))=\e(w_i\T^l+(1+w)\T^l_i)\left(\n_{X_l}X_\al-h_{l\al}\calV^{\e}
\right),
$$
$$
\frac{\pa}{\pa y^a} X_\al(  y,(1+w(y,z))
\TT+\e^{-1}\Phi(\e y))=\e\delta_{ab}\left(\n_{X_b}X_\al-h_{b\al}\calV^{\e} \right)+\e\left( w_a\T^l
 +\Phi_a^l \right)\left(\n_{X_l}X_\al-h_{l\al}\calV^{\e} \right).
$$
\begin{proposition}\label{pr:sff}
The following expansions hold
\begin{eqnarray}\label{eq:nyaZa}\nonumber
 \la  N, \frac{\pa}{\pa y^a} Z_{{a}} \ra & = &-\e\Gamma_{a}^a (\TT )
+\e\Tm h_{aa}-w_{aa}-\e \,\la \Phi_{ a  a},\TT  \ra -\e \,\la R
(\Phi , E_a) \, E_a , \TT \ra
 \\
&+& \e\, \Gamma_{a}^c (\TT) \,\Gamma_{c}^a (\Phi)-2\e \Tm w_a h(\TT)^a
+\frac{\e}{\mu^2}w_l\left( \Gamma_{a}^a (\TT_l )-h_{aa}\Theta_l^{m+1}   \right)\\
&+&\mathcal{O}(\e^2)+ \e^2 \, L (w, \Phi)  + \e  Q(w, \Phi);\nonumber
\end{eqnarray}

\begin{eqnarray}
\la N, \frac{\pa}{\pa z^j} Z_{j} \ra & = & \mu^2 (1+w) -w_{jj}-\e \Tm \la h(\TT),\TT \ra
w_{jj}-2\e\Theta_j^{n+1}\la h(\TT),\TT \ra  w_j\nonumber\\
&+& \e (1+2 w) \left(\Tm \la h(\TT_j),\TT_j \ra
-2 \Theta_j^{n+1} \la h(\TT),\TT_j \ra -\Tm \la h(\TT),\TT_{jj} \ra \right)\nonumber\\[3mm]
&+&\frac{\e}{\mu^2}w_k \left(2\Tm \la h(\TT_k),\TT_{ii} \ra +2\Theta_i^{n+1}
\la h(\TT_k),\TT_{i} \ra + \Theta_k^{n+1}
\la h(\TT_i),\TT_{i} \ra  \right)\nonumber
\\[3mm]
& + & \frac{2}{3} \, \e \,  \la  R(\Phi , \TT_j ) \, \TT,
\TT_j \ra -\frac{\e}{3} \, \la R(\Phi, \TT) \, \TT , \TT_{jj}\ra
+2\e w_c \left( \Theta_j^{n+1}h(\TT_j)^c+\Tm h(\TT_{jj})^c   \right)\nonumber\\[3mm]
& + &2 \e \la \Phi_{\bar c},\TT \ra \left( \Theta_j^{n+1}h(\TT_j)^c+\Tm h(\TT_{jj})^c \right)
+ \e \Tm h(\TT)^c\left( \la \Phi_{ c},\TT_{jj}\ra +\mu^2\la \Phi_{c},\TT \ra  \right) \\[3mm]
&+& \e \Tm h(\TT)^c\left(w_c \la \TT,\TT_{jj}\ra +\mu^2w_c  \right)
-\frac16 \mu^2 \la R(\Phi, \TT) \, \Phi , \TT\ra -\frac13\la R(\Phi, \TT) \, \Phi , \TT_{jj}\ra \nonumber\\[3mm]
&-& \frac 12 \mu^2 w_c^2+\frac 12\mu^2 |\la\Phi_{ c},\TT \ra|^2
  -\frac 12 w_k^2+ 2 w_j^2 +\la \Phi_{ c},\TT_{jj} \ra w_c
   +\la \Phi_{c},\TT \ra  \la \Phi_{ c},\TT_{jj}\ra \nonumber\\[3mm]
&+& (1+2w) \al^k\la\T_{jj},\T_k  \ra +\mathcal{O}(\e^2)+ \e^2 \, L (w, \Phi)  + \e  Q(w, \Phi);\nonumber
\end{eqnarray}
\begin{eqnarray*}
\la  N, \frac{\pa}{\pa y^a} Z_{{b}} \ra & = &
-\Gamma_{a}^b (\TT ) +\e\Tm h_{ab} -w_{ab}+{\cal O}(\e^2)+ \e L (w, \Phi)+ Q(w, \Phi))
\quad a \neq b; \\[3mm]
\la  N, \frac{\pa}{\pa {y}_a} Z_{j} \ra & = & \e \Theta_j^{n+1} h(\TT)^a+\e\Tm h(\TT_j)^a -w_{aj}
+{\cal O}(\e^2)+ \e L (w, \Phi)+ Q(w, \Phi); \\[3mm]
\la  N, \frac{\pa}{\pa z_i} Z_{j} \ra & = & -w_{ij}-\e \Theta_i^{n+1}
\la h(\TT),\TT_j \ra  -\e \Theta_j^{n+1} \la h(\TT),\TT_i \ra +\e \Tm
\la h(\TT_i),\TT_{j} \ra\\[3mm]
&-&\e \Tm
\la h(\TT),\TT_{ij} \ra +\al^k \la \T_{ij},\T_k \ra + {\cal O} (\e^2) +
 \e L(w, \Phi)+ Q(w, \Phi), \quad i \neq j.
\end{eqnarray*}
\end{proposition}
{\bf Proof~:} The proof is similar in spirit to the one of Proposition 3.3 in \cite{mmp}.
So we will be sketchy here referring to the aforementioned paper for more details.
We have that
\begin{eqnarray*}
  \frac{\pa }{\pa y^a}Z_a &=& \e \left( \n_{X_a}X_a-h_{aa}{\mathcal{V}^{\e}} \right)
  +w_{aa}\Upsilon +2\Tm w_a D_a {\mathcal{V}^{\e}}+\e \Phi_{{a}{a}}^l X_l
  +\Tm D_aD_a {\mathcal{V}^{\e}}\\
  &+& \left( \mathcal{O}(\e^2)+\e^2L(w,\Phi)+\e Q(w,\Phi)  \right)X_\al+\left(\mathcal{O}(\e^2)
  +\e^2 L(w,\Phi)+\e Q(w,\Phi)
   \right)\calV^\e
\end{eqnarray*}
and for $a\ne b$
\begin{eqnarray*}
 \frac{\pa }{\pa y^a}Z_b &=& \e \left( \n_{X_b}X_a-h_{ab}{\mathcal{V}^{\e}} \right)
  +w_{ab}\Upsilon\\[3mm] &+& \left( \mathcal{O}(\e^2)+\e L(w,\Phi)+ Q(w,\Phi)  \right)X_\al
  +\left(\mathcal{O}(\e^2)+\e^2 L(w,\Phi)+\e Q(w,\Phi)
   \right)\calV^\e;
\end{eqnarray*}

\begin{eqnarray*}
  \frac{\pa }{\pa z^i}Z_i &=& w_{ii}\Upsilon+2\,w_i \Upsilon_i
  +2\e \T^l\T_i^s w_i\left(\n_{X_s}X_l-h_{sl}{\mathcal{V}^{\e}} \right)
  +2\Tm D_i{\mathcal{V}^{\e}}\,w_i+(1+w)\Upsilon_{ii}\\[3mm]
  &+& (1+w)\left( 2\Theta_i^{n+1}D_i{\mathcal{V}^{\e}}+\Tm D_iD_i\mathcal{V }  \right)
  +\e(1+2w)\T_i^l\T_i^s \left(\n_{X_s}X_l-h_{sl}{\mathcal{V}^{\e}} \right)\\[3mm]
  &+& \left( \mathcal{O}(\e^2)+\e^2L(w,\Phi)+\e Q(w,\Phi)  \right)X_\al+\left(\mathcal{O}(\e^2)
  +\e^2 L(w,\Phi)+\e Q(w,\Phi)
   \right)\calV^\e;
\end{eqnarray*}
and for $i\ne j$
\begin{eqnarray*}
  \frac{\pa }{\pa z^i}Z_j &=& w_{ij}\Upsilon+w_i \Upsilon_j+w_j \Upsilon_i
  +\Theta_i^{n+1}D_j{\mathcal{V}^{\e}}
  +\Theta_j^{n+1}D_i{\mathcal{V}^{\e}}+(1+w)\Upsilon_{ij}+
  \e \T_i^l\T_j^s \left(\n_{X_s}X_l-h_{sl}{\mathcal{V}^{\e}} \right)\\
  &+&\Tm D_iD_j\mathcal{V }^{\e}+ \left( \mathcal{O}(\e^2)+\e L(w,\Phi)
  + Q(w,\Phi)  \right)X_\al+\left(\mathcal{O}(\e^2)+\e^2 L(w,\Phi)+\e Q(w,\Phi)
   \right)\calV^\e.
\end{eqnarray*}
Finally
\begin{eqnarray*}
\frac{\pa }{\pa y^a}Z_j=\frac{\pa }{\pa z_j}Z_a &=&
 \e \T_j^s\left( \n_{X_s}X_a-h_{as}{\mathcal{V}^{\e}} \right)
  +w_{aj}\Upsilon +w_a\Upsilon_j+\Theta^{n+1}_j D_a {\mathcal{V}^{\e}}\\
  &+& \left( \mathcal{O}(\e^2)+\e^2L(w,\Phi)+\e Q(w,\Phi)  \right)X_\al
  +\left(\mathcal{O}(\e^2)+\e^2 L(w,\Phi)+\e Q(w,\Phi)
   \right)\calV^\e.
  \end{eqnarray*}

 \

\noindent Recalling the expansions, see Lemma 2.1 in \cite{mmp}.
\begin{equation}
\begin{array}{rllll}
\nabla_{X_i} \, X_j  & = &   ( {\cal O}
(\e) + L(w, \Phi)+ Q(w, \Phi))  \, X_\gamma, \\[3mm]
\nabla_{X_a} \, X_i  & = & - \Gamma^b_{a} (E_i) \, X_b + ( {\cal O}
(\e) + L (w, \Phi)+ Q(w, \Phi))  \, X_\gamma,\label{eq:3-7}
\end{array}
\end{equation}
We will also need the following expansion which follows from the
result of Lemma 2.2 in \cite{mmp} (with obvious modifications).
\begin{eqnarray}\nonumber
\nabla_{X_a} \, X_b & = &  \Gamma_{a}^b (E_j) \, X_j -  \la R( \e \, \TT +
 \Phi , E_a) \, E_j , E_b \ra \, X_j \\[3mm]
& + & \frac{1}{2} \, \left( \la R (E_a, E_b)\, (\e \, \TT +
\Phi) , E_j\ra - \Gamma^c_{a} (\e \, \TT + \Phi) \,\Gamma_{c}^b
(E_j) - \Gamma^b_c(\e \, \TT + \Phi) \,
\Gamma_{a}^c (E_j) \right)  \, X_j   \\[3mm]
& + & ( {\cal O} (\e) + L(w, \Phi) + Q(w, \Phi) ) \, X_c + (
{\cal O} (\e^2 ) +  \e \, L(w, \Phi) + Q (w, \Phi) ) \,
X_j.\nonumber
\label{eq:3-2t}
\end{eqnarray}
These implies in particular
\begin{eqnarray*}
  \la \Upsilon,\n_{X_a}X_a\ra &=& \T^l\Gamma_a^a(E_i)\left( \de_{li}+2\e\Tm h_{li} \right)
  -\e\la R(\TT,E_a)\TT,E_a \ra -\la R(\TT,E_a)\Phi,E_a \ra\\
  &-& \e \Gamma_a^c(\TT)\Gamma_c^a(\TT)-\Gamma_a^c(\TT)\Gamma_c^a(\Phi)+\mathcal{O}(\e^2)
  +\e L(w,\Phi)+Q(w,\Phi).
\end{eqnarray*}
On the other hand we have that
\begin{equation*}
 D_a D_a{\mathcal{V}^{\e}}=\e w_{aa}h(\TT)^\al X_\al +\left(\mathcal{O}(\e^2)+\e^2 L(w,\Phi)+\e Q(w,\Phi)
   \right)X_\be+\left(\mathcal{O}(\e^2)+\e^2 L(w,\Phi)+\e Q(w,\Phi)
   \right)\calV^\e,
\end{equation*}
which implies
\begin{equation}
\la D_a D_a{\mathcal{V}^{\e}}, \Upsilon\ra =\e w_{aa}\la h(\TT),\TT\ra
+\mathcal{O}(\e^2)+\e^2 L(w,\Phi)+\e Q(w,\Phi).
\end{equation}
Using these together with \eqref{eq:bc}, \eqref{eq:alj} and Lemma \ref{l:ovg},
 the first estimate follows at once.
For the other estimates one can proceed similarly.\hfill $\Box$

\subsubsection{The mean curvature of perturbed tubes}\label{sss:mc}

Collecting the estimates of the last subsection we obtain the
expansion of the mean curvature of the hypersurface $S_\e(w ,
\Phi)$. In the coordinate system defined in the previous sections,
we get
\begin{equation}
\begin{array}{rlllllll}
 m \, H (w, \Phi) & = & n - \e \, \Gamma^a_a (\TT) +
\e\,\Tm\, h_{aa}+\e\,\Tm \left[ (n+3)\la h(\TT),\TT \ra -h_{jj}\right]
+ {\cal O} (\e^2) \notag \\[3mm]
  &-& \big(  \Delta_{K_\e} w + \Delta_{S^{n}} w + n
  w \big) -  \e  \bigg( \la\, \Delta_K \Phi +  R(\Phi, E_a) \, E_a ,
  \TT \,\ra - \Gamma_a^c (\Phi) \, \Gamma_c^a (\TT)\bigg)  \\[3mm]
  &-& \e \Tm \la h(\TT),\TT \ra \, \Delta_{S^n}w -2\e (n+3)\,\Tm \la h(\TT),\nabla_{S^n}w \ra
 + 2\e \Tm\,\nabla^2_{S^n}w:h  \\[3mm]
 &-& \e \left(\la h(\TT),\TT \ra + h_{jj}
 +h_{aa}\right)\la \nabla_{S^n}w,E_{n+1}\ra
 -(1+3n)\e \Tm h(\TT)^a\,w_a\\[3mm]
 &-& 2\e \Tm h(\nabla_{S^n}w_a)^a+\e\, \Gamma^a_a (\nabla_{S^n}w)
  -2\e \nabla^2_{K_\e}w:\Gamma(\TT)
 +2\e \Tm h_{aa} w_{aa} \\[3mm]
 &-& (3n+1)\e\,\Tm h(\TT)^a \la \Phi_{\bar{a}},\TT \ra+\e\,\Tm h(\Phi_{\bar{a}})^a
  +2\e\Tm h:\Gamma(\Phi)\\[3mm]
 &+& n\,w^{2}+\frac{2-n}{2}\left| \nabla_{S^n}w\right|^2+2\,w\,\Delta_{S^n}w
 -\frac n2(w_a+\la \Phi_{\bar{a}},\TT \ra)^2          \\[3mm]
&-& \la\Phi_{\bar a},\nabla_{S^n}w_a \ra-2 \nabla^2_{K_\e}w:\Gamma(\Phi) +
\frac{n+2}{6}\la R(\Phi,\TT)\Phi\,,\,\TT \ra-\frac 13 \la R(\Phi,E_i)\Phi\,,\,E_i \ra\\[3mm]
& + & \mathcal{O}(\e^2)+ \e^2 \,  L(w, \Phi) +\e\,Q(w, \Phi). \notag \label{eq:mc}
\end{array}
\end{equation}
Here we have used the formulas in Lemma \ref{l:formulas}, the fact that
$$
\Delta_{S^n}=\frac{1}{\mu^2}\left( \Delta_{\R^n}+(2-n)\T^i \pa_{i}  \right),
$$
and the notation
$A:B=A_{st}B_{st}$ for two linear operators $A$ and $B$.
Where summation over repeated indices is understood.
We first define the following operators appearing in the above expansion
\begin{eqnarray}\label{eq:defL1}\nonumber
\mathcal{L}^1(w):&=&-\la h(\TT),\TT \ra \, \Delta_{S^n}w -2 (n+3)\,\Tm \la h(\TT),\nabla_{S^n}w \ra
 + 2 \Tm\,\nabla^2_{S^n}w:h  \\[3mm]
 &-&  \left(\la h(\TT),\TT \ra + h_{jj}+h_{aa}\right)\la \nabla_{S^n}w,E_{n+1}\ra\\[3mm]
 &-&\e(1+3n)\Tm \la h(\TT),\nabla_{K} w\ra+\e\Tm  h(\nabla_{S^n}w_{\bar{a}})^a
 -2  \e^2\nabla^2_{K}w:\Gamma(\TT)+2\e^2 \Tm h_{aa} w_{aa}
 ,\nonumber
\end{eqnarray}
%
%

\begin{equation}\label{eq:J1}
\mathcal{J}^1\Phi:=-(3n+1)\,\Tm h(\TT)^a \la \Phi_{\bar{a}},\TT \ra+\,\Tm h(\Phi_{\ov{a}})^a
  +2\Tm h:\Gamma(\Phi),
\end{equation}
and the quadratic term
\begin{eqnarray}\label{eq:Q}\nonumber
\mathcal{Q}^1(w,\Phi):&=&n\,w^{(2)}+\frac{2-n}{2}\left| \nabla_{S^n}w\right|^2
+2\,w\,\Delta_{S^n}w-\frac n2(\e w_{\bar{a}}+\la \Phi_{\bar{a}},\TT \ra)^2
-\e\la\Phi_{\bar{a}},\nabla_{S^n}w_{\bar{a}} \ra\\
&-&2\e^2 \nabla^2_{K}w:\Gamma(\Phi) +
\frac{n+2}{6}\la R(\Phi,\TT)\Phi\,,\,\TT \ra-\frac 13 \la R(\Phi,E_i)\Phi\,,\,E_i \ra.
\end{eqnarray}

\medskip

Next, we define
\[
{\cal L}_\e : = \e^2 \, \Delta_K + \Delta_{S^{n}} + n,\qquad\qquad {\cal L}_0 : =  \Delta_{S^{n}} + n
\]
and the Jacobi operator about $K$ in $(\pa\Om,\bar{g})$, see \S~\ref{ss:Jac-K}
\[
{\mathfrak J} : = \Delta^{\perp} - {\cal R}^{\perp} + {\mathcal B}.
\]
 Recall that (see \S~\ref{ss:stgp}) the outer unit normal to the boundary
 of $\pa S^n_+$ in $S^n_+$ is   $\eta=-E_{n+1}$,
%
$$
\frac{\pa w}{\pa \eta}=-\la\n_{S^n_+}\,w,E_{n+1}\ra.
$$
Using these definitions, we obtain the following result~:
\begin{proposition}\label{pr:4.1}
Assume that $K$ is a minimal submanifold, then the mean curvature
of ${ S}_\e (w, \Phi)$ can be expanded as
\begin{equation}
\begin{array}{rlllllll}
 m \, H (w, \Phi) & = &  n + \e\,\Tm\, h_{aa}+\e\,\Tm
\left[ (n+3)\la h(\TT),\TT \ra -h_{jj}\right]+{\cal O} (\e^2)\\[3mm]
&-& {\cal L}_\e \, w-  \e \,  \la {\mathfrak J} \, \Phi, \TT\ra
+\e \mathcal{L}^1 w+ \e \mathcal{J}^1(\Phi)+\mathcal{Q}^1(w,\Phi)\\[3mm]
&+&\e^2 \, L(w, \Phi) +\e\, Q(w, \Phi). \notag
\end{array}
\label{eq:mcb}
\end{equation}
where $\mathcal{L}^1$ is defined in \eqref{eq:defL1}, $\mathcal{J}^1$
is given in \eqref{eq:J1}
while $\mathcal{Q}^1$ is a quadratic term  defined in \eqref{eq:Q}.
Moreover, the orthogonality condition is equivalent to the
following boundary condition on the function
$w$:
\begin{equation}\label{eq:orthog}
\frac{\pa w}{\pa \eta}=w\,\frac{\pa w}{\pa \eta}+\bar{{\cal O} }(\e^2)
+\e^2\, \bar{L}(w,\Phi)+\e\, \bar{Q}(w,\Phi)\quad \hbox{ on } \pa (SNK)_+.
\end{equation}
\end{proposition} {\bf Proof :} The expression of the mean curvature can be obtained
 rather easily taking into account the above definitions (with obvious modifications)
 and the minimality of $K$ which implies
\[
\Gamma_a^a =0.
\]
\hfill $\Box$

\medskip
With these notations finding $w$ and $\Phi$ such that the equation $
m \, H= n$ and $\la N,\calV^\e\ra=0$ hold is equivalent to solve
\begin{equation}\label{ndf}\left\{
\begin{array}{rllll}
{\cal L}_\e \, w + \e \,  \la {\mathfrak J} \, \Phi, \TT\ra &=& \e\,\Tm\, h_{aa}+\e\,\Tm
\left[ (n+3)\la h(\TT),\TT \ra -h_{jj}\right]+ {\cal O} (\e^2)\\[3mm]
& +&  \e\, \mathcal{J}^1(\Phi)+\e \mathcal{L}^1 w+\mathcal{Q}^1(w,\Phi)+
 \e^2 \, L(w, \Phi) +\e\, Q(w, \Phi)\quad\hbox{ in }  (SNK)_+,\\[3mm]
\displaystyle \frac{\pa w}{\pa \eta}&=&\displaystyle w\,\frac{\pa w}{\pa \eta}
+\bar{{\cal O}} (\e^2)+\e^2\, \bar{L}(w,\Phi)+\e\, \bar{Q}(w,\Phi)\qquad\qquad
\qquad\quad\hbox{ on } \pa (SNK)_+.
\end{array}\right.
\end{equation}
%
%


%
%
\section{Adjusting the tube $\bar{S}_\e
(K_\e)$}\label{s:appsol}

In this section we annihilate the error terms ($\calO(\e)$) appearing in
(\ref{ndf}) at any given order. The non-degeneracy of the  submanifold $K$
will play a crucial
 role in such a construction.
 We denote by  $\Pi$ the $L^2$  projection on the subspace spanned by the $\Theta^i$, $i=1,\cdots ,n$.\\
We set  $$\hat{w}^{(r)}=\sum_{d=1}^r\e^d w^{(d)} \qquad\textrm{ and }\qquad\hat{\Phi}^{r}=\sum_{d=1}^{r-1} \e^d \Phi^{(d)}.$$ \\
\textbf{Construction of $w^{(1)}$:} We first want to kill the term
$\calO(\e)$. This is equivalent to having
\begin{eqnarray}\label{eq:eqw2}
\left\{
  \begin{array}{ll}
   m\,H(\hat{w}^{(r)},\hat{\Phi}^{(r)})=n+\calO(\e^2),\quad\textrm{ in }\quad S_\e
(\hat{w}^{(r)},\hat{\Phi}^{(r)}),\\[3mm]
  \la N,{\mathcal{V}^{\e}}\ra=\bar{\calO}(\e^{2})\quad\textrm{ on }\quad\pa S_\e
(\hat{w}^{(r)},\hat{\Phi}^{(r)}).
  \end{array}
\right.
\end{eqnarray}
This gives the following equation in $w^{(1)}$
\begin{equation}\label{eq:eqw1}
\begin{array}{rllll}
\mathcal{L}_0 w^{(1)}&=&\,\Tm\, h_{aa}+\Tm
\left[ (n+3)\la h(\TT),\TT \ra -h_{jj}\right]\quad\textrm{ in }(SNK)_+;\\
\displaystyle \frac{\pa w^{(1)}}{\pa \eta}&=&0\quad\textrm{ on }\pa(SNK)_+.\displaystyle
\end{array}
\end{equation}
By the result from \S~\ref{ss:stgp} (with $\gamma=\frac{\pi}{2}$) and  Fredholm
alternative theorem, the solvability of \eqref{eq:eqw1}
is possible provided
$$
\int_{S^n_+}\left(\Tm\, h_{aa}+\Tm
\left[ (n+3)\la h(\TT),\TT \ra -h_{jj}\right]\right)\,\Theta^i\,d\te=0\qquad\textrm{ for all } i=1,\cdots, n
$$
which is the case by oddness, here $d\te$ denotes the volume element on $S^n_+$.\\
 Notice that the variable $\bar{y}$ is being considered as a parameter so that $w^{(1)}$
 is as smooth as the right hand side in this  variable.

\

\noindent  \textbf{Constructing $w^{(2)}$ :} We turn now to the term of
order $\e^2$. We have
\begin{eqnarray}\label{eq:eqw2}
\left\{
  \begin{array}{ll}
   mH(\hat{w}^{(r)},\hat{\Phi}^{(r)})=n+\calO(\e^3),\quad\textrm{ in }\quad S_\e
(\hat{w}^{(r)},\hat{\Phi}^{(r)}),\\
   \la N,{\mathcal{V}^{\e}}\ra=\bar{\calO}(\e^{3})\quad\textrm{ on }\quad\pa S_\e
(\hat{w}^{(r)},\hat{\Phi}^{(r)}).
  \end{array}
\right.
\end{eqnarray}
Since the terms involving  $\Phi$ in $\mathcal{Q}^1(\e w^{(1)},\e \Phi^{(1)})$ are $\e^3 L(\Phi^{(1)})$
 and $ Q(\hat{\Phi}^{(r)},\hat{\Phi}^{(r)})$, (\ref{eq:eqw2})  yields a system in $w^{(2)}$ and $\Phi^{(1)}$ given by
\begin{equation}\label{eq:solvw2}
\begin{array}{rllll}
\mathcal{L}_0 w^{(2)}&=&\la \mathfrak{J}\Phi^{(1)},\TT \ra+\mathcal{O}(1)+
\mathcal{L}^1 w^{(1)}+ \mathcal{J}^1(\Phi^{(1)})+Q(\hat{\Phi}^{(r)},\hat{\Phi}^{(r)})\quad\textrm{ in }(SNK)_+\\
\displaystyle \frac{\pa w^{(2)}}{\pa \eta}&=&\bar{\mathcal{O}}(1)\quad\textrm{ on }\pa(SNK)_+.
\end{array}
\end{equation}
 Note that $\Pi\,\mathcal{J}^1=0$ and $\Pi\,Q(\Phi^{(1)},\Phi^{(1)})=0$  so (\ref{eq:solvw2})
 is solvable if and only if
$$
\int_{S^n_+}\la \mathfrak{J}\Phi^{(1)},\TT \ra\,\Theta^i\,d\te+ \int_{S^n_+}\left( \mathcal{O}(1)
+  \mathcal{L}^1 w^{(1)}\right)\,\Theta^i\,d\te
+   \oint_{\pa S^n_+} \bar{\mathcal{O}}(1)\,\Theta^i\,d\bar{\te }=0\quad \textrm{for all } i=1\cdots n,
$$
where $d\te$ and  $d\bar{\te}$ are the volume elements on $S^n_+$
and $\pa S^n_+$ respectively. This gives an equation on $\Phi^{(1)}$
which can be solved using the non degeneracy of the submanifold $K$,
once this is done, the solvability on $w^{(2)}$ follows at once.

\

\noindent \textbf{Constructing $w^{(r)}$:} We want to construct an
approximate solution as accurate as possible, and to do so we will use
an iterative scheme. Suppose the couple $(w^{(r-1)},\Phi^{(r-2)})$ is
already determined. To find $(w^{(r)},\Phi^{(r-1)})$, it suffices to
check that when we project on the Kernel of $\mathcal{L}_0$, the
operator involving $\Phi^{(r-1)}$ should be  only the Jacobi operator
$\mathfrak{J}$. This is the case since the only term that can bring
 $\Phi^{(r-1)}$ at this iteration step is
$\mathcal{Q}^1_{r-1}(w,\Phi)$ which gives only
 terms of the form $\e^2\Phi$ and  $Q(\hat{\Phi}^{(r)},\hat{\Phi}^{(r)})$
 moreover $\Pi\,\mathcal{J}^1_{r-1}(\Phi^{(r-1)})=\Pi\,Q (\hat{\Phi}^{(r)},\hat{\Phi}^{(r)})=0$.\\
The index $r$ appearing in the
 linear and quadratic terms means that they depend on the iteration step
 while the operator $\mathcal{J}^1_r$  keep its same properties because it
 is influenced only by the even quadratic terms in $Q(\hat{\Phi}^{(r)}+\Phi,\hat{\Phi}^{(r)}+\Phi)$
  appearing in $Q^1(\hat{w}^{(r)}+w,\hat{\Phi}^{(r)}+\Phi)$.\\
 %
 \\
By induction, in the same argument, for every $r\in \N$, we can find
$(w^{(d)},\Phi^{(d)})$, $d=1,\cdots, r$ smooth such that
\begin{equation}\label{le:lklk}
\hat{w}^{(r)}=\sum_d^r\e^d w^{(d)}=\calO(\e)\quad\textrm{ and }\quad\hat{\Phi}^{(r)}=\sum_d^{r-1}\e^d \Phi^{(d)}=\calO(\e)
\end{equation}
and that
$$
m\,H(\hat{w}^{(r)},\hat{\Phi}^{(r)})= n+\calO(\e^{r+1})\quad \textrm{ in } \quad S_\e
(\hat{w}^{(r)},\hat{\Phi}^{(r)}),\qquad \la N,{\mathcal{V}^{\e}}\ra
=\bar{\calO}(\e^{r+2})\quad\textrm{ on }\quad\pa S_\e(\hat{w}^{(r)},\hat{\Phi}^{(r)}).
$$

\begin{remark}
 Notice that as in \cite{mm} we omitted the terms involving
derivatives with respect to $\bar y$ of the function $w$ (by considering $\mathcal{L}_0$
 instead of $\mathcal{L}_\e$), this is due to
fact that since $w$ is slow dependent on $y_a$,  when differentiating with respect
to $y_{\bar a}$ we pick up an $\e$ at each differentiation, this gives us smaller
terms. However, when applying elliptic regularity theorems we might loose two
derivatives at each iteration. This indeed is not a problem since one needs
just a finite number of iterations. We refer the reader to \cite{mm}, where
a more explanation is given.
\end{remark}
We are left to find  $w$ and $\Phi$ such that
\begin{equation}\label{eq:tosovle}
\begin{array}{ccc}
 m\,H(\hat{w}^{(r)}+w,\hat{\Phi}^{r}+\Phi)&=&n \qquad\textrm{ in }\quad S_\e
(\hat{w}^{(r)}+w,\hat{\Phi}^{r}+\Phi),\\[3mm]\qquad \la N,{\mathcal{V}^{\e}}\ra&=&0\qquad\textrm{on }\quad\pa S_\e
(\hat{w}^{(r)}+w,\hat{\Phi}^{r}+\Phi).
\end{array}
\end{equation}
We define the linearized mean curvature operator about $S_\e (\hat{w}^{r}, \hat{\Phi}^{r})$
$$
\mathbb{L}_{\e,r}(w,\Phi)=\frac{1}{\e}\left({\cal L}_\e \, w
+\e\mathcal{ L}^{1}_r(w)\right) +  \,  \langle {\mathfrak J}\Phi,\TT\rangle
+ \mathcal{J}^1_r({\Phi})+\e L_r(w,\Phi).
$$
The index $r$ appearing in the constant,
 linear and quadratic terms means that they depend on the iteration step but keep there properties.\\
We Notice that $\mathbb{L}_{\e,r}$ is not precisely the usual
Jacobi operator because we are parametrizing this hypersurface as a graph
over $S_\e (\hat{w}^{r}, \hat{\Phi}^{r})$ using the vector field $-
\Upsilon$ rather than the unit normal $N$.

\medskip
Using Remark  \ref{rem:linkJac} ($\gamma=\frac{\pi}{2}$),  suppose that
$\Sigma = S_\e (\hat{w}^{r}, \hat{\Phi}^{r})$ and
$\hat{N} = - \Upsilon$. From (\ref{le:lklk}) and
Proposition~\ref{pr:normal} we have
\[
\la N, - \Upsilon\ra = 1 + {\cal O} (\e^2).
\]
Furthermore, from Proposition~\ref{pr:metric} and (\ref{le:lklk}),
 the volume forms of the tubes
$S_\e (\hat{w}^{r}, \hat{\Phi}^{r})$ and $(SN K)_+$ are related by
\[
dvol_{S_\e (\hat{w}^{r}, \hat{\Phi}^{r})} =   (1+ {\cal
O} (\e)) \,  dvol_{(SN K)_+}.
\]
We define $\delta_{\e,r} >0$ by
\begin{equation}
\la N, -\Upsilon\ra \, dvol_{S_\e (\hat{w}^{r}, \hat{\Phi}^{r})} =
   \delta_{\e ,r} \, dvol_{(SN K)_+}. \label{eq:NNN}
\end{equation}
Multiplying by $\delta_{\e,r}$, the system (\ref{eq:tosovle}) will change
 the terms $\calL^1_r$, $L_r$, $\bar{L}_r$,  the constant and  quadratic terms
  will keep there properties and there will be a new linear operator $\bar{\calL}^{1}_r(w)$
  on the boundary. We keep the same notations for these terms and call $\mathbb{L}_{\e,r}$
  the new selfadjoint operator $\delta_{\e,r}\,\mathbb{L}_{\e,r}$ with respect to
  the standard $L^2(SNK)_+$-inner product.\\


%
Now since $\bar{L}_r(w,\Phi)$ and $\bar{\calL}^{1}_r(w)$ involves only terms of the
form $w$, ${\pa_{z^i} w}$, the boundary
conditions can be changed to
$$
(1+\calO(\e))\,\frac{\pa w}{\pa \eta}=\e\, \hat{L}_r(w,\Phi)+\hat{\calL}^{1}_r(w)
+\bar{\mathcal{O}}_r(\e^{r+1})  + \frac{1}{\e}w\,\frac{\pa w}{\pa \eta}
+\, \bar{Q}_r(w,\Phi)\quad \hbox{ on } \pa (SNK)_+
$$                                                                                                                     %
with $\hat{L}_r(w,\Phi) $ and  $\hat{\calL}^{1}_r(w)$ contain no angular
derivatives in $w$. Now by the trace theorem we can extend
$\hat{L}_r(w,\Phi)$, $\hat{\calL}^{1}_r(w)$  and $\bar{\mathcal{O}}_r(\e^{r+1})$
 in $(SNK)_+$  and this will just add some terms in ${L}_r(w,\Phi)$, ${\calL}^{1}_r(w)$
   and ${\mathcal{O}_r}(\e^{r})$ respectively  which will maintain there properties.
We conclude that there is no loss of generality  when replacing the solvability of
 (\ref{eq:tosovle})  with the following equation.
\begin{equation}\label{eq:solvph}
\begin{array}{rllll}
 \mathbb{L}_{\e,r}(w,\Phi)&=&\frac{1}{\e} Q_r(w,\Phi)+\mathcal{O}_r(\e^{r})\quad \hbox{ in } (SNK)_+,\\[3mm]
 \displaystyle \frac{\pa w}{\pa \eta}&=&\displaystyle \frac{1}{\e} \bar{Q}_r(w,\Phi)\,\,
 \quad\quad\qquad \hbox{ on } \pa (SNK)_+.
\end{array}
\end{equation}
We will try to  invert  the linear operator on the left hand side and this will lead us to
 study the spectrum of the operator by selfadjointness.
\section{Spectral analysis}\label{s:span}
\textbf{Function space:} Fix $\frac{1}{2}>s>0$. For any $v \in
L^2(SNK)_+$, set
\[
 \la\Phi,\TT\ra:=\Pi \, v , \qquad \e^{-1+2s} \, w:=\Pi^\perp v ,
\]
so that
\begin{equation}\label{eq:dec-v}
v = \e^{1-2s} \, w + \la\Phi,\TT\ra.
\end{equation}
It will be understood that $\Phi^i$ for $i=1,\cdots, n$ are the components of $\Pi\,v$ on $NK$.
Conversely if $(w,\Phi)\in \Pi^{\perp}\,L^2(SNK)_+\times L^{2}(K,NK)$ is given,
we associate to it $v$ as in (\ref{eq:dec-v}).

\medskip

Later  we will often  decompose
\begin{equation}
w = w_0 + w_1\label{eq:0mode}
\end{equation}
where $w_0$ is a function on $K$ and  $w_1$ has zero mean value with respect to the angular integrals.

\medskip
The volume element of $(SN K)_{+}=S^n_+\times K$ will be denoted by $d\theta\, d\bar{y}$.\\
As it will be apparent later, we will be considering the following weighted Hilbert subspaces of $L^2(SNK)_+$
$$
L^2_{\e}:=\left\{v=\e^{1-2s} \, w + \la\Phi,\TT\ra\in L^2(SNK)_+\quad:
\quad \e^{-2s}\,\int_{(SN K)_{+}}  |w|^2 \, d\te\, d\bar{y}
+ \int_{K} |\Phi|^2 \, d\bar{y}<\infty\right\}
$$
with corresponding norm
$$
\|v\|^2_{L^2_{\e}}:=\e^{-2s}\,\int_{(SN K)_+}  |w|^2 \, d\te \,d\bar{y} + \int_{K} |\Phi|^2 \, d\bar{y}.
$$
We also define
\begin{eqnarray*}
H^1_{\e}:=\left\{v\in L^2_{\e}\quad:\quad \e^{-2s}\,\int_{SN K} (\e^2 \, |\nabla_K w|^2 +
|\nabla_{S^{n}} w|^2  + |w|^2) \, d\te\,d\bar{y}
+ \int_{K} ( |\nabla_K \Phi|^2 + |\Phi|^2) \, d\bar{y}<\infty\right\}
\end{eqnarray*}
with  corresponding norm
$$
\|v\|_{H^1_{\e}}^2:= \e^{-2s}\,\int_{SN K} (\e^2 \, |\nabla_K w|^2 +
|\nabla_{S^{n}} w|^2  + |w|^2) \, d\te\,d\bar{y}
+ \int_{K} ( |\nabla_K \Phi|^2 + |\Phi|^2) \, d\bar{y}.
$$
Let $ |S^{n}_+|$ denote the volume of $S^{n}_+$.
Notice that
$$
\int_{S^{n}_+}(\Theta^i)^2\, d\te=\frac{|S^{n}_+|}{n+1}\qquad \textrm{for all }i=1\cdots n.
$$
We define $\varrho_n:=\frac{|S^{n}_+|}{n+1}$.
\medskip
With these definitions in mind we redefine ${\mathbb L}_{\e ,r}$ by duality as follows
\begin{eqnarray*}
 && \int_{(SN K)_+} v  \, {\mathbb L}_{\e
,r} \, v' \, d\te\,d\bar{y}: =  \\[3mm]
& - &\e^{-2s}\int_{(SN K)_{+}} \e^2w' \, \, \Delta_K w
\,d\te\,d\bar{y}+ \e^{-2s}\int_{(SN K)_{+}}( \nabla_{S^{n}_+}\,
 w\,\nabla_{S^{n}_+}\, w' - n \, w\,w') \, d\te\,d\bar{y}\\[3mm]
& + & \varrho_n\, \int_{K}  \la\mathfrak
J \Phi, \Phi'\ra \, d\bar{y}+ \int_{(SNK)_+} (  \mathcal{J}^1_r({\Phi})+\mathcal{L}^1_r(w)
+ \e L_r (w,\Phi)) \, (\e^{1-2s}\, w'+ \la\Phi', \TT\ra) \, \,d\te\,d\bar{y}.
\end{eqnarray*}
We associate to ${\mathbb L}_{\e ,r}$ its quadratic bilinear form
\[
{\cal C}_{\e ,r} (v,v') : = \int_{(SN K)_+} v  \, {\mathbb L}_{\e
,r} \, v' \, d\te\,d\bar{y},
\]
and the associated quadratic form ${\cal Q}_{\e ,r} (v) : = {\cal
C}_{\e ,r} (v,v)$.\\
As mentioned in the first section, following \cite{mm}, we want to find
the values of $\e$ for which the operator ${\mathbb L}_{\e,r}$ is invertible.
 By selfadjointness this leads to find the values of $\e$ for which the
 eigenvalues of the form ${\cal Q}_{\e ,r}$ are bounded away from zero.
 Such techniques requires first that  our form should be very close to a
  model one that we can characterize its  spectrum (just the small eigenvalues).
  Secondly, to understand the behavior of small eigenvalues seeing as ``set" valued functions in $\e$.
   We will  estimate the Morse index of ${\cal Q}_{\e ,r}$ and prove  the monotonicity of its
    small eigenvalues. The former  can be done using Weyl's asymptotic formula and the latter
     can be obtained by applying a result by  Kato. We shall do this in the remaining of this
      section.

\

\noindent We define  the model form, by duality, as
\begin{eqnarray*}
{\cal C}_0 (v,v')  &: =&  \displaystyle - \e^{-2s}\int_{(SN K)_{+}}\e^2 w' \,
\, \Delta_K w \,d\te\,d\bar{y}
+ \e^{-2s}\int_{(SN K)_{+}}( \nabla_{S^{n}_+}
\, w\,\nabla_{S^{n}_+}\, w' - n \, w\,w') \,
d\te\,d\bar{y}\\
&+&\varrho_n\, \int_{K}  \la\mathfrak
J \Phi, \Phi'\ra \, d\bar{y}
\end{eqnarray*}
and the associated quadratic form ${\cal Q}_0  (v) : = {\mathcal C}_0
(v,v)$.

\begin{proposition}\label{p:model}
There exists a constant $c >0$ (independent of $r$) such that
\begin{equation}
\left|{\cal C}_{\e ,r} (v, v') - {\cal C}_0 (v,v') \right| \leq c
\, \e^s \, \|v\|_{H^1_\e} \, \|v'\|_{H^1_\e}.\label{eq:5.55}
\end{equation}
\end{proposition}
{\bf Proof~:} First of all we notice that in  $\mathcal{
L}^{1}_r(w)$  may appear expressions of the forms
$w$, $\e\del_{\ov{y}^a} w$, $\e^2 \,
\del_{\ov{y}^a} \del_{\ov{y}^b}w$, $\del_{z^j} w$, $\del_{z^j} \del_{z^{j'}}w$.
 Nevertheless after integrating by parts and using H\"older inequality there holds
$$
 \left| \int_{(SNK)_{+}}\e^{1-2s} w'\,\mathcal{ L}^{1}_r(w) \,d\te\,d\bar{y}\right|\leq\e c\| v\|_{H^1_{\e}}\,\|
 v'\|_{H^1_{\e}},
$$
  and by definition of the $H^1_{\e}$ norm
 \begin{eqnarray*}
 \left| \int_{(SNK)_+}\la\Phi',\TT\ra\,\mathcal{ L}^{1}_r(w)
 \,d\te\,d\bar{y}\right| &\leq&c \e^{s} \| \e^{1-2s} w\|_{H^1_{\e}}\, \|
\Phi'\|_{L^2(K,NK)}\\[3mm]
&\leq& c\e^{s}\| v\|_{H^1_{\e}}\,\| v'\|_{H^1_{\e}}.
 \end{eqnarray*}
Furthermore  $\Pi \mathcal{J}^1({\Phi})=0$.
Now it is clear that  even if $ \mathcal{J}^1_r({\Phi}) +  L_r (w,\Phi)$
involves terms of the form $w$, $\e\del_{\ov{y}^a} w$, $\e  \,
\del_{\ov{y}^a} \del_{\ov{y}^b}w$, $\del_{z^j} w$, $\del_{z^j} \del_{z^{j'}}
w$ and also $\Phi^j$, $ \del_{\ov{y}^a} \Phi^j$ and $\del_{\ov{y}^a}
\del_{\ov{y}^b} \, \Phi^j$, in any case  after integration by parts and using H\"older inequality,
$$
\left|\int_{(SNK)_+} ( \e^{-1}\, \mathcal{J}^1_r({\Phi})+  L_r (w,\Phi)) \, (\e^{1-2s}\, w'
+ \la\Phi', \TT\ra) \, \,d\te\,d\bar{y}\right|\leq c\| v\|_{H^1_{\e}}\,\| v'\|_{H^1_{\e}}.
$$
The result follows at once. \hfill $\Box$

\

\noindent \textbf{The Morse index of ${\cal Q}_{\e ,r}$:} Define the
two quadratic forms
\[
{\cal Q}^{\pm} (v) : =  {\cal Q}_0 (v) \pm  \gamma \, \e^s \, \| v
\|_{H^1_\e}^2 .
\]
From (\ref{eq:5.55}), if $\gamma
>0$ is sufficiently large and $\e$ small enough, then
\[
{\cal Q}^- \leq {\cal Q}_{\e ,r} \leq {\cal Q}^+,
\]
so that the index of ${\cal Q}_{\e ,r}$ is bounded by those of ${\cal Q}^+$ and ${\cal Q}^-$.

\medskip

Given any function $w$ defined on $(SNK)_+$, we set
\[
D^\pm_0(w) : = (1 \pm \gamma \, \e^s) \, \int_{K} \e^2 \,
|\nabla_K w|^2 \,d\bar{y} - (n \mp \gamma \, \e^s) \, \int_{K}
|w|^2 \, d\bar{y},
\]
\[
D^\pm_1(w) : = (1 \pm \gamma \, \e^s ) \, \int_{(SN K)_+} (\e^2 \,
|\nabla_K w|^2 + |\nabla_{S^{n}_+}  w|^2 ) \, d\te\,d\bar{y}- (n \mp
\gamma \, \e^s) \, \int_{(SN K)_+} |w|^2 \, d\te\,d\bar{y},
\]
and finally,
\[
D^\pm(\Phi) : = - ( 1\pm \gamma \, \e^s) \, \int_K \la{\mathfrak J}
\, \Phi, \Phi\ra \, d\bar{y}.
\]
With these definitions in mind, we have
\[
{\cal Q}^\pm (v) = (n+1)\varrho_n\, \e^{-2s}\, D^\pm_0 (w_0) + \e^{-2s}\,D^\pm_1(w_1)+
\varrho_n \,  D^\pm (\Phi) ,
\]
if we decompose $v = \e^{1-2s} \, w+ \la\Phi, \TT\ra$ and further
decompose $w = w_0+w_1$ as usual. Following Section 6.3 in \cite{mmp}
it is easy to see that if $( 1\pm \gamma \, \e^s)>0$ then the index of $D^\pm$ is
 the index of $K$. Moreover the index of $D^\pm_1$ is equal to zero if   $2 \,
(n+1) \, (1 - \gamma\, \e^s) - (n  + \gamma \, \e^s) >0$ because
\[
\Pi \, w_1 = 0\qquad \mbox{and} \qquad  \int_{S^{n}_+} w_1  \,
d\te = 0
\]
hence
\[
\int_{S^{n}_+} |\nabla_{S^{n}_+} w_1|^2 \, d\te \geq 2\, (n+1)
\, \int_{S^{n}_+} |w_1|^2 \, d\te.
\]
This shows that the asymptotic behavior of the index of ${\cal Q}_{\e ,r}$
should be determined by $D_0^\pm$. It is the case since   its index is given by
$$
\sharp\{j\quad:\quad (1\pm\gamma\e^s)\lambda_j<(n\mp\gamma\e^s)\},
$$
 where $\lambda_j$ are the eigenvalues of $-\e^2\Delta_K$ counted with multiplicities.
Now using Weyl's formula one obtain its index,
$$
\mbox{Ind} \, D_0^\pm \sim c_K \, \left( \frac{n}{\e^2}
\right)^{\frac{k}{2}}.
$$
Collecting these estimates, one obtains the following
\begin{lemma}\label{lem:index}
The Morse index of ${\cal Q}_{\e ,r}$ is asymptotic to $c\e^{-k}$ when $\e$
 tends to zero, where $c$ depends only on $m$ and $K$.
\end{lemma}

\noindent \textbf{Approximate eigenfunctions: } In order to apply
Kato's theorem \cite {Kato} we need to characterize the
eigenfunctions (eigenspaces) corresponding to small eigenvalues. We
prove
\begin{lemma} Let
$\sigma$ be an eigenvalue of $\,{\mathbb L}_{\e ,r}$ and $v = \e^{1-2s} \, w + \la\Phi, \TT \ra$
a corresponding eigenfunction and  $\e^{1-2s}\,w_0 =\int_{S^n_+}v\,d\te$ is the decomposition from (\ref{eq:0mode}).
 There
exist constants $c, c_0 >0$ such that if $|\sigma| \leq c_0$, then
\[
\|v - \e^{1-2s}\,w_0\|_{H^1_\e}^2 \leq c\,\e^s\,\|v\|_{H^1_\e}^2,
\]
for all $\e >0$ small enough.
 \label{le:local}
\end{lemma}
{\bf Proof:} For any $v' = \e^{1-2s} \, w' + \la\Phi' , \Theta \ra$, we have
\[
\begin{array}{llll}
\calC_{\e ,r} (v,v') & = & \displaystyle \sigma \int_{SN K} (
\e^{2-4s} w \, w'  + \la\Phi,\Theta\ra \la\Phi',\Theta\ra ) \, d\te\,d\bar{y}\\[3mm]
& = & \displaystyle \sigma \int_{SN K} \e^{2-4s} w \, w' \, d\te\,d\bar{y}
+ \sigma \, \varrho_n \, \int_K \la\Phi,\Phi'\ra\,d\bar{y}.
\end{array}
\]
In addition, (\ref{eq:5.55}) gives
\begin{equation}
\begin{array}{lll}
\displaystyle \left| \int_{SN K} \e^{-2s}(\e^2 \, \nabla_K w \, \nabla_K
w' +\nabla_{S^{n}_+} w  \nabla_{S^{n}_+} w' - (n +\sigma \, \e^{2-4s}) \, w \, w') \, d\te\,d\bar{y}\right. \\[3mm]
\displaystyle\qquad\qquad\qquad\qquad+\left.\varrho_n  \, \int_{K} (\la\mathfrak J \Phi, \Phi'\ra -
\sigma \, \la\Phi, \Phi'\ra ) \,  d\bar{y}\right| \leq c \, \e^s \,
\|v\|_{H^1_\e} \, \|v'\|_{H^1_\e}.
\end{array}
\label{eq:fdfd}
\end{equation}

\noindent {\bf Step 1 : }
Let $\Phi' =0$
and $w' = w_1$ to get
\[
\left| \int_{SN K} \e^{-2s}(\e^2 \, |\nabla_K w_1|^2 + |\nabla_{S^{n}_+}
w_1|^2 - (n-\sigma \, \e^{2-4s}) \, |w_1|^2  ) \,  d\te\,d\bar{y}\right|
\leq c \, \e^s \, \|v\|_{H^1_\e}\,\| \e^{1-2s} \, w_1 \|_{H^1_\e}.
\]
However, since
\[
\Pi \, w_1 = 0\qquad \mbox{and} \qquad  \int_{S^{n}_+} w_1  \,
d\te = 0,
\]
we have
\[
\int_{S^{n}_+} |\nabla_{S^{n}_+} w_1|^2 \, dvol_{S^{n}_+} \geq 2\, (n+1)
\, \int_{S^{n}_+} |w_1|^2 \, d\te,
\]
hence
\[
\left| \int_{SN K} \e^{-2s}(\e^2 \, |\nabla_K w_1|^2 + \frac{1}{2} \,
|\nabla_{S^{n}_+} w_1|^2 + (1 - |\sigma| \, \e^{2-4s}) \, |w_1|^2 ) \,d\te\,
d\bar{y}\right| \leq c \,\e^s\,\|v\|_{H^1_\e}^2.
\]
This implies that
\[
\| \e^{1-2s}\,w_1 \|_{H^1_\e}^2  \leq c \, \e^s  \,
\|v\|_{H^1_\e}^2 ,
\]
for all $\e \in (0,1)$, provided $|\sigma | \leq 1/2$.\\

\noindent {\bf Step 2: } Now let $w'=0$ and $\Phi' =\Phi^+$ (resp.
$\Phi' = \Phi^-$) in (\ref{eq:fdfd}), where $\Phi^+$ (resp.
$\Phi^-$) is the $L^2$ projection of $\Phi$ over the space of
eigenfunctions of $\mathfrak J$ associated to positive (resp.
negative) eigenvalues. This yields
\[
\left| \int_{K} ( \la \mathfrak J \Phi, \Phi^\pm\ra - \sigma \, \la\Phi
, \Phi^\pm\ra ) \,  d\bar{y}\right| \leq c \, \e^s  \, \|v\|_{H^1_\e}
\, \|\la\Phi^\pm , \TT\ra\|_{H^1_\e}.
\]
Since ${\mathfrak J}$ is invertible, there exists $c_1 >0$ such that
\[
c_1 \, \| \la\Phi^\pm , \TT\ra\|^2_{H^1_\e} \leq \left|
\int_{K} \la \mathfrak J \Phi, \Phi^\pm\ra \,  d\bar{y}\right|.
\]
Hence
\[
( c_1 - |\sigma |) \, \|\la\Phi^\pm , \TT\ra\|_{H^1_\e}^2 \leq
c \, \e^s \, \|v\|_{H^1_\e}^2.
\]
This conclude the proof with $c_0:=\min\{1/2,c_1/2\}$.

 \hfill $\Box$
 \begin{remark}\label{rm:cont-nH}
 If $v$ is in an eigenfunction corresponding to an eigenvalue given by the above lemma,  then it satisfies
\begin{equation*}
\begin{array}{rlllll}
\displaystyle \left| \int_{SN K}\e^{-2s} (\e^2 \, |\nabla_K w|^2 +
|\nabla_{S^{n}_+} w|^2 - (n + \sigma\e^{2-4s} ) \, |w|^2)  \,
d\te\,d\bar{y}\qquad \qquad\right. \\[3mm]
\displaystyle \hfill +\left.\varrho_n \, \int_{K} ( \,\la
\mathfrak J \Phi, \Phi\ra - \sigma\, \la\Phi, \Phi\ra ) \, d\bar{y}
\right| \leq c \, \e^s \, \|v\|_{H^1_\e}^2,
\end{array}
\end{equation*}
and
\begin{equation}\label{eq:us-var}
 \left| \int_{SN K}\e^{-2s} (\e^2 \, |\nabla_K w|^2 +
|\nabla_{S^{n}_+} w|^2 - n   \, |w|^2)  \,
 d\te\,d\bar{y}\right| \leq c \, \e^s \, \|v\|_{H^1_\e}^2.
\end{equation}
Notice  that $\nabla_{S^{n}_+} w = \nabla_{S^{n}_+} w_1$ if $w$ is
decomposed as $w=w_0+ w_1$ one has
\begin{equation*}
\left| \int_{SN K} \e^{-2s}(\e^2 \, |\nabla_K w|^2 - n \,
|w|^2)  \,  d\te\,d\bar{y}\right| \leq  c \, \e^s \, \|v\|_{H^1_\e}^2,
\end{equation*}
so that
$$
\e^{-2s}\int_{SN K} \e^2 \, |\nabla_K w|^2\,  d\te\,d\bar{y}\leq  c \, \e^s \,
 \|v\|_{H^1_\e}^2+n\e^{-2s}\int_{SN K} |w|^2  \, d\te\,d\bar{y}.
$$
In particular we have
$$
\|v\|_{H^1_\e}\leq c \|v\|_{L^2_\e}.
$$
\end{remark}

\noindent \textbf{Variation of small eigenvalues with respect to
$\e$:} To understand the behavior of small eigenvalues of the
symmetric quadratic form $\mathcal{Q}_{\e,r}$, we need to apply a result by Kato, see \cite{Kato}.
 Considering the eigenvalues $\sigma(\e)$ as differentiable multivalued function  in $\e$. The result states that
\begin{equation}\label{eq:Kato-res}
\del_\e  \sigma \in \left\{ \int_{SN K} v \, (\del_\e {\mathbb
L}_{\e, r} ) \,  v \, d\te\,d\bar{y} \qquad  : \qquad {\mathbb
L}_{\e ,r} v = \sigma \, v, \qquad \| v \|_{L^2} =1 \right\}.
\end{equation}
In order to obtain some informations  about the spectral gaps of the
linearized operator when the parameter $\e$ is small, one can look
at its small eigenvalues as differentiable function on $\e$,
differentiate them with respect to $\e$ and estimate their
derivatives. This is indeed given in the following Lemma.
\begin{lemma}
There exist constants $c_1,c >0$ such that, if $\sigma$ is an
eigenvalue of ${\mathbb L}_{\e ,r}$ with $|\sigma| < c_1$, then
\[
\e \, \del_\e \sigma \geq 2 \, n - c \, \e^s ,
\]
provided $\e$ is small enough. \label{le:evol}
\end{lemma}
{\bf Proof~:}
We have just to  provide bounds for the set on the right  of (\ref{eq:Kato-res}) using the above remark.
\medskip
Assume that ${\mathbb L}_{\e ,r} v = \sigma \, v$, but rather than
normalizing the function $v$ by $\|v \|_{L^2} =1$, assume instead
that $\|v\|_{L_\e^2} =1$. In order to compute $\del_\e {\mathbb
L}_{\e ,r}$, recall that
\[
w = \e^{-1+2s} \, \Pi^\perp v \qquad \mbox{and that} \qquad
\la{\mathfrak J}\Phi,\TT\ra = \Pi \, v,
\]
so we can write
\begin{eqnarray*}
{\mathbb L}_{\e ,r} \, v &=&  - \e^{2s}\Delta_K  \, (\Pi^\perp  \, v) +
\frac{1}{\e^{2-2s}} \, \calL_0 \, (\Pi^\perp  \, v) + \Pi  \, v
+ \frac{1}{\e^{1-2s}}\,\mathcal{L}^1_r\,(\Pi^\perp  \, v)\\
& &\qquad\qquad\qquad +\mathcal{J}^1\,({\mathfrak J}^{-1}_r\Pi\, v)
+ \e \, {L}_r( {\e}^{-1+2s} \, \Pi^\perp v  ,{\mathfrak J}^{-1} \Pi \, v).
\end{eqnarray*}
Since $\Pi$ and $\Pi^\perp$ are independent of $\e$, we have
\begin{eqnarray*}
\del_\e {\mathbb L}_{\e ,r} v& =& - 2s\e^{-1+2s}\Delta_K  \, (\Pi^\perp  \, v)
 + {(-2+2s)}{\e^{-3+2s}} \, \calL_0 \,
(\Pi^\perp v)  +{(-1+2s)}{\e^{-2+2s}} \,\mathcal{L}^1_r\,(\Pi^\perp  \, v) \\[3mm]
&+& \tilde L_r ( \e^{-1+2s} \, \Pi^\perp v,  {\mathfrak J}^{-1}\Pi\, v ),
\end{eqnarray*}
where the operator $\tilde L_r$  varies from line to line but
satisfies the usual assumptions. This now gives
\begin{eqnarray*}
\left| \int_{SN K} v \, (\del_\e {\mathbb L}_{\e ,r} )  \, v \,
{d\te\,d\bar{y}} \right.&-&\left. {2}\e^{-1-2s}\int_{SN K}\e^2|\nabla_K w|^2\,{d\te\,d\bar{y}}\right.\\[3mm]
 &+&\left.\frac{(2-2s)}{\e} \, \e^{-2s}\int_{SN K} (\e^2|\nabla_K w|^2+|\nabla_{S^{n}_+} w|^2-
n \, |w|^2)  \, {d\te\,d\bar{y}}\right|\\[3mm]
 &\leq& c \, \| v \|^2_{H^1_\e}+\left|\frac{1-2s}{\e}\,\int_{SN K}\,
 \la\Phi,\TT\ra \mathcal{L}^1_r\,(w)\,{d\te\,d\bar{y}}\right|\\[3mm]
&\leq& \frac{c}{\e^{1-s}} \, \| v \|^2_{H^1_\e}
 .
\end{eqnarray*}
Consequently if $v$ is an eigenfunction of ${\mathbb L}_{\e, r}$ with
corresponding eigenvalue $|\sigma|\leq c_0$, where $c_0$ is given in
the previous lemma, by the inequality \eqref{eq:us-var}, see the above remark, we have
\begin{equation}
\left| \int_{SN K} v \, (\del_\e {\mathbb L}_{\e ,r} )  \, v \,
 {d\te\,d\bar{y}}- {2}\e^{-1-2s}\int_{SN K}\e^2|\nabla_K w|^2\,{d\te\,d\bar{y}}
\right|\leq \frac{c}{\e^{1-s}} \, \| v \|^2_{H^1_\e}.
 \label{eq:cc}
\end{equation}
Using again from the above remark, one gets 
\[
\e^{-1-2s}\int_{SN K}\e^2|\nabla_K w|^2\,{d\te\,d\bar{y}}\leq c\,\e^{-1+s}
 \, \|v\|_{H^1_\e}^2+n\,\e^{-1-2s}\,\int_{SN K} |w|^2  \, {d\te\,d\bar{y}}.
\]
If we normalize $v$ by $\|v\|_{L^2_\e} =1$
 then inserting this into (\ref{eq:cc}) we get
\begin{equation}
\left| \int_{SN K} v \, (\del_\e {\mathbb L}_{\e ,r} )  \, v  \,
 {d\te\,d\bar{y}} -  \frac{2}{\e} \, n \right|\leq  \frac{c}{\e^{1-s}}\label{frfr}
\end{equation}
for all eigenfunction $v$ such that ${\mathbb L}_{\e , r} v =
\sigma \, v$ which is normalized by $\| v \|_{L_\e^2} =1$.

\medskip

This already implies that $\del_\e \sigma >0$ for $\e$ small
enough. But observing that we always have $||v||_{L^2} \leq
\|v\|_{L_\e^2}$, we conclude that
\[
\inf_{ \stackrel{{\mathbb L}_\e v = \sigma \, v}{\|v\|_{L^2} = 1}}
\int_{SN K} v \, (\del_\e {\mathbb L}_\e ) \, v  \, {d\te\,d\bar{y}}
\quad \geq \inf_{\stackrel{{\mathbb L}_\e v = \sigma \,
v}{\|v\|_{L_{\e}^2} = 1}} \int_{SN K} v \, (\del_\e {\mathbb
L}_\e) \, v \, {d\te\,d\bar{y}},
\]
and (\ref{frfr}) implies that
\[
\del_\e \sigma \geq  \frac{2}{\e} \, n - \frac{c}{\e^{1-s}}.
\]
This completes the proof of the result. \hfill $\Box$
\section{Proof of Theorem \ref{th:existence}}
Using  Lemma \ref{lem:index}   and Lemma \ref{le:evol}, reasoning as for
 the proof of Lemma 6.3 in \cite{mmp} we can find a sequence of open
 interval $I_i$, $i\in\N$ such that the smallest eigenvalue of
 ${\mathbb L}_{\e ,r}$ is bounded away from zero for any $\e\in\cup_i I_i$. More precisely we have
\medskip

\begin{lemma}
Fix any $q \geq 2$. Then there exists a sequence of disjoint
nonempty open intervals $I_i = (\e_i^-, \e_i^+)$,
$\e_i^\pm \rightarrow 0$ and a constant $c_q > 0$ such that
when $\e \in I^q : =  \cup_i I_i$, the operator ${\mathbb
L}_{\e ,r}$ is invertible and
\[
({\mathbb L}_{\e ,r})^{-1} : L^2_\e \longrightarrow L^2_\e
,
\]
has norm bounded by $c_q \,  \e^{-k-q+1}$, uniformly in $\e \in
I$. Furthermore, $I^q :=  \cup_i I_i$ satisfies
\[
\left|{\mathcal H}^1 ( (0, \e) \cap I^q) -\e \right|\leq c \,
\e^{q}, \qquad \e \searrow 0.
\]
\label{le:sg}
\end{lemma}
For $p\in\N$ and $0<\al<1$, we denote by $\mathcal{C}^{p,\al}$ the usual H\"older spaces  on the closure of $(SNK)_+$.
\begin{lemma}\label{lem:invert}
Let $f\in \mathcal{C}^{0,\al}$  and $v$ satisfy
$$
{\mathbb L}_{\e ,r}\,v=f.
$$
Then there exit a constant $c>0$ (independent of $\e$ but depend on $r$)
and  $R>0$ depending only on $q$, $\al$, $s$ and $k$ such that
$$
\|v\|_{\mathcal{C}^{2,\al}}\leq c\,\e^{-R}\,\|f\|_{\mathcal{C}^{0,\al}}
$$
for any $\e\in I^q$.
\end{lemma}
{\bf Proof~:}
Fix $q\geq 2$. Observe that by definition of the weighted norm of $L^2_\e$, from Lemma \ref{le:sg} we have
$$
\|v\|_{L^2}\leq c_q \,  \e^{-k-q+1-s}\,\|f\|_{L^2}.
$$
By standard elliptic regularity theory, there exists $c>0$ (depending on $r$)
 such that  the following H\"older estimate holds
$$
\e^{2+\al}\|v\|_{\mathcal{C}^{2,\al}}\leq c\,\e^2\,\|f\|_{\mathcal{C}^{0,\al}}+c\,\e^{-\frac{k}{2}}\,\|v\|_{L^2}.
$$
From these last two inequalities, we can  choose $R>\frac{3k}{2}+q+\al+1+s$.
\hfill $\Box$\\
\\
\\
We end the proof of the main theorem by finding a fixed point for the mapping
$$
T_{\e,r}(v):=-({\mathbb L}_{\e ,r})^{-1}\left\{\mathcal{O}_r(\e^r)+\mathcal{N}_{\e,r}(v)\right\},
$$
where
$$
\begin{array}{ccc}
\displaystyle \int_{(SNK)_+}\mathcal{N}_{\e,r}(v)\,v'\,d\te\,d\bar{y}
&:=&\displaystyle\int_{(SNK)_+}\e^{-1}\,Q_r(\e^{-1+2s}\,\Pi^{\perp}\,v,\Pi\,v)\,v'\,d\te\,d\bar{y}\\[3mm]
 &&+\displaystyle\oint_{\pa(SNK)_+}\e^{-1}\,\bar{Q}_r(\e^{-1+2s}
 \,\Pi^{\perp}\,v,\Pi\,v)\,v'\,d\bar{\te}\,d\bar{y}.
\end{array}
$$
Since by definition, $Q_r$ and $\bar{Q}_r$  are (at least) quadratic we have
\begin{eqnarray*}
\|\mathcal{N}_{\e,r}(v)\|_{\mathcal{C}^{0,\al}}&=& \e^{-2+2s}\,
O(\|v\|_{\mathcal{C}^{2,\al}})\,\|v\|^2_{\mathcal{C}^{2,\al}};\\[3mm]
\|\mathcal{N}_{\e,r}(v_1)-\mathcal{N}_{\e,r}(v_2)\|_{\mathcal{C}^{0,\al}}
&=& \e^{-2+2s}\,O(\|v_1\|_{\mathcal{C}^{2,\al}},  \|v_2\|_{\mathcal{C}^{2,\al}})
\|v_1-v_2\|_{\mathcal{C}^{2,\al}}.
\end{eqnarray*}
%
  Now we fix $ r>2\, R+2-2\,s$. By Lemma \ref{lem:invert} and the above inequalities,
   for every $\e\in I^q$, $T_{\e,r}(v)$  maps the ball
  $$
  \{v\in\mathcal{C}^{2,\al}\quad:\quad \|v\|_{\mathcal{C}^{2,\al}}\leq C\,\e^{r+1-R}\}
  $$
into itself moreover it is a contraction. Therefore it has a unique
 fixed point $v=\e^{{1-2s}}\,w+\la\Phi,\TT\ra$ in the ball yielding
\begin{equation*}
\begin{array}{cccc}
 m\,H(\hat{w}^{(r)}+w,\hat{\Phi}^{r}+\Phi)&=&n& \quad\textrm{ in }\quad S_\e
(\hat{w}^{(r)}+w,\hat{\Phi}^{r}+\Phi)\subset\Om_\e,\\[3mm]
\qquad \la N,{\mathcal{V}^{\e}}\ra&=&0&\quad\textrm{ on }\quad\pa S_\e
(\hat{w}^{(r)}+w,\hat{\Phi}^{r}+\Phi)\subset\pa\Om_\e.
\end{array}
\end{equation*}
If $\e\in I^q$ is sufficiently small then  rescaling back, the tube $\,\e\, S_\e
(\hat{w}^{(r)}+w,\hat{\Phi}^{r}+\Phi)$,  is an embedded hypersurface of $\Om$
with constant mean curvature equal to $\frac{n}{m}\e^{-1}$ and
intersecting the boundary of $\Om$ perpendicularly along its boundary.
\\
\begin{remark}\label{rem:ex-cap}\textbf{Existence of stationary Capillary hypersurfaces}.\\
\\
Letting $\gamma\in(0,\pi)$ be an angle, recall from \S~\ref{ss:fc} that
 $(y^1,y^{2}\,\dots\,,y^k)\in\R^k$ (resp. $(z^1,z^{2}\,\dots\,,z^n)\in B^n_{r(\gamma)}$)
  are the local coordinate variables on $K_\e$ (resp. on $S^n(\gamma)$), where
    $r(\gamma):=\frac{1-\cos{\gamma}}{1+\cos(\gamma)}$  (see \S~\ref{ss:stgp}) and
$$
\T(\gamma):=\bfp{\Big|_{B^n_{r(\gamma)}}}-\cos(\gamma)\,E_{n+1}
$$
parametrize the spherical cap $S^n(\gamma)$ which intersect the horizontal
plane $\R^m$ with angle $\gamma$.\\
As in the case where $\gamma=\frac{\pi}{2}$, we can use the same class of deformations letting
$\Phi:K\to NK_\e$ and $w:B^n_\gamma\times K_\e\to \R$, consider
$$
S_\gamma\;:(y,z)\mapsto y\times \e^{-1}\Phi(\e y)\,+ \,(1+w(y,z))\,\T(\gamma).
$$
The nearby surfaces  around $K_\e$  which make an angle
almost equal to $\gamma$ with $\pa\Om_\e$ can be parametrized (locally) by
$$
G_{\gamma}(y,z):(y,z)\longrightarrow S_{\gamma}(y,z)\longrightarrow F^\e(S_\gamma(y,z)),
$$
namely
\[
G_\gamma(y,z) :=  F^{\e} \, \left( { y},\frac{1}{\e}\Phi(\e y)
+(1+w(y,z))\tilde{\T}(\gamma),(1+w(y,z))\T^{n+1}(\gamma)\right).
\]
Notice that $\Tm(\gamma)\Big|_{\pa B^n_{r(\gamma)}}=0$, so there holds
$$
G_\gamma(y,z)\Big|_{\pa B^n_{r(\gamma)}}\in\pa\Om_\e\qquad\textrm{ for any } y
$$
 The image of this map will be called
$S_\e^{\gamma}(w,\Phi)$. \\
Observe that the hypersurfaces close to $S_\e^{\gamma}(0,0)$
 are parametrized using the vectorfield $-\Upsilon(\gamma) = \Theta^j (\gamma)\, X_j +\T^{n+1}(\gamma)\mathcal{V}^{\e}$
  rather than the normal $\Xi:= \bfp^j \, X_j +\bfp^{n+1} \mathcal{V}^{\e}$
because it  is more  reasonable if we want  the boundary of
 $S_\e^{\gamma}(w,\Phi)$ to be on $\pa\Om_\e$ without imposing
 simultaneously a Neumann and Dirichlet boundary condition on $w$. Suppose $Z_j(\gamma), Z_a(\gamma)$ span the tangent space of $S_\e^{\gamma}(w,\Phi)$
as in  \S~\ref{sss:nvf},  we can obtain the  normal fields $N(\gamma)$ by finding $\al^j(\gamma)$
and $\be^{a}(\gamma)$ so that
 $$
 N(\gamma)=-\Xi+\al^j(\gamma) Z_j(\gamma)+\be^a(\gamma) Z_a(\gamma).
 $$
As we did so far, the mean curvature at every point  of $S_\e^{\gamma}(w,\Phi)$ can be easily obtained
\begin{equation}
\begin{array}{lllllll}
 m \, H (w, \Phi) & = & n - \e\bigg( \, \Gamma^a_a (\tilde{\bfp}) +\bfp^{n+1}\, h_{aa}
 +\bfp^{n+1} \left[ 3\la h(\tilde{\bfp}),\tilde{\bfp} \ra -h_{jj}\right]
 +n\,\T^{n+1}(\gamma) \la h(\tilde{\bfp}),\tilde{\bfp} \ra\bigg)+{\cal O} (\e^2)
\notag \\[3mm]
  &-& \bigg( \e^2 \Delta_{K}\left(\la\T(\gamma),\bfp\ra w \right)
  + \Delta_{S^{n}} \left(\la\T(\gamma),\bfp\ra w \right)
   + n   \left(\la\T(\gamma),\bfp\ra w \right)\bigg) \\[3mm]
  &-&  \e  \bigg( \la\, \Delta_K \Phi +  R(\Phi, E_a) \, E_a\, ,\,
  \tilde{\bfp} \,\ra - \Gamma_a^c (\Phi) \, \Gamma_c^a (\tilde{\bfp})\bigg)  \\[3mm]
 &-&\e\bigg( (3n+1)\,\Tm(\gamma) h(\tilde{\bfp})^a \la \Phi_{\bar{a}},\tilde{\bfp} \ra
 +\,\bfp^{n+1} h(\Phi_{\bar{a}})^a   +2\bfp^{n+1} h:\Gamma(\Phi)\bigg)\\[3mm]
 &-& \frac n2(\e w_{\bar{a}}+\la \Phi_{\bar{a}},\tilde{\bfp} \ra)^2
  -  \la\Phi_{\bar a},\e\nabla_{S^n}w_{\bar{a}} \ra
  -2 \e^2\nabla^2_{K}w:\Gamma(\Phi)  \\[3mm]
&+ &\frac{n+2}{6}\la R(\Phi,\tilde{\bfp})\Phi\,,\,\tilde{\bfp} \ra
-\frac 13 \la R(\Phi,E_i)\Phi\,,\,E_i \ra\\[3mm]
& +&  \e\, L(w)+ \e^2 \,  L(w, \Phi)+Q(w)  +\e\,Q(w, \Phi). \notag \label{eq:mcgam}
\end{array}
\end{equation}
Moreover (recall that $\calV^{\e}$ is the interior normal of $\pa\Om_\e$)
 using the fact that $\Tm(\gamma)\Big|_{\pa B^n_{r(\gamma)}}=0$,   the equation
  $\la-\calV^{\e},N\ra= \cos(\gamma)$ is equivalent to
\begin{eqnarray*}
\la\T(\gamma),\bfp\ra(1-w)\frac{\pa w}{\pa \eta(\gamma)}
 =\bar{{\cal O} }(\e^2)+\e^2\, \bar{L}(w,\Phi)
 +\bar{Q}^1(w,\Phi)+\e\, \bar{Q}(w,\Phi)\qquad \hbox{ on } \pa S^n(\gamma)\times K,
\end{eqnarray*}
which is again equivalent to
\begin{eqnarray*}
\frac{\pa (\la\T(\gamma),\bfp\ra w)}{\pa \eta(\gamma)}
 &=&w\frac{\pa \la\T(\gamma),\bfp\ra }{\pa \eta(\gamma)}
 +\bar{{\cal O} }(\e^2)+\e^2\, \bar{L}(w,\Phi)+\bar{Q}^1(w,\Phi)+\bar{Q}(w)\\[3mm]
 &&+\e\, \bar{Q}(w,\Phi)\qquad \qquad \qquad \qquad \qquad \qquad \qquad \qquad \hbox{ on } \pa S^n(\gamma)\times K\\[3mm]
&=&w\,\cot(\gamma)+\bar{{\cal O} }(\e^2)+\e^2\,
\bar{L}(w,\Phi)+\bar{Q}^1(w,\Phi) \\[3mm]
&&+\bar{Q}(w)+\e\, \bar{Q}(w,\Phi)\,\,\,\,\qquad \qquad \qquad \qquad \qquad \qquad \hbox{ on } \pa
S^n(\gamma)\times K,
\end{eqnarray*}
where
$$
\bar{Q}^1(w,\Phi):=\cot(\gamma)\left( \e w_{\bar{a}} \la \Phi_{\bar{a}},\tilde{\bfp} \ra
+ \la \Phi_{\bar{a}},\tilde{\bfp} \ra \la \Phi_{\bar{a}},\tilde{\bfp} \ra
   -\frac{1}{3} \, \la R (\Phi, \tilde{\bfp}) \, \Phi , \tilde{\bfp} \ra    \right).
   $$
Using the result from \S~\ref{ss:stgp} and from \S~\ref{s:appsol}, one can
adjust the tube to  $S_\e^{\gamma}(\hat{w}^{(r)},\hat{\Phi}^{(r)})$ accurately.
Moreover with the decomposition of the  functions $v=\e^{1-2s}\, w
+\la\Phi,\tilde{\bfp}\ra\in L^2(S^n(\gamma)\times K)$  as in
(\ref{eq:dec-v}) we conclude that the spectral analysis of the linearized
mean curvature operator over $S_\e^{\gamma}(\hat{w}^{(r)},\hat{\Phi}^{(r)})$ carried
out as we obtain in Section \ref{s:span} in the new weighted Hilbert subspaces
 of $L^2(S^n(\gamma)\times K)$
\begin{eqnarray*}
L^2_{\e,\gamma}&:=&\bigg\{v=\e^{1-2s} \, w +
\la\Phi,\tilde{\bfp}\ra\in L^2(S^n(\gamma)\times K)\quad:\\
&&\qquad \quad \e^{-2s}\,\int_{S^n(\gamma)\times K} \la\T(\gamma),
\bfp\ra |w|^2 \, d\te(\gamma)\, d\bar{y} + \int_{K} |\Phi|^2 \,
 d\bar{y}<\infty\bigg\}
\end{eqnarray*}
%
%
%
\begin{eqnarray*}
&&\bigg\{v\in L^2_{\e,\gamma}\quad:\quad
\e^{-2s}\,\int_{S^n(\gamma)\times K}
 \la\T(\gamma),\bfp\ra(\e^2 \, |\nabla_K w|^2 +
|\nabla_{S^{n}} w|^2  + |w|^2) \, d\te(\gamma)\,d\bar{y}\\
&&\qquad \qquad\qquad\qquad\qquad + \int_{K} ( |\nabla_K \Phi|^2 +
|\Phi|^2) \, d\bar{y}<\infty\bigg\}.
\end{eqnarray*}
Under the usual assumptions on $K$, if $\e\in I^q$ is sufficiently small then
 rescaling back, the tube $\,\e\, S_\e
(\hat{w}^{(r)}+w,\hat{\Phi}^{r}+\Phi)\,$,  is an embedded  hypersurface of $\Om$ with
constant mean curvature $\frac{n}{m}\e^{-1}$ and intersecting $\pa\Om$ with and angle $\gamma$.
 This yields  a set of stationary Capillary hypersurfaces in $\Om$ with constant
 ``contact angle" $\gamma$ and condensing to the  submanifold $K$.

\end{remark}

\

\begin{center}
{\bf Acknowledgments}
\end{center}

\noindent The authors wish to thank Professor Andrea Malchiodi for
useful discussions and Professor Frank Pacard for pointing out the
generalization to  capillary problems. They are supported by
M.U.R.S.T within the PRIN 2006 Variational Methods and Nonlinear
Differential Equations and by GNAMPA within the project 2007
 Geometric Evolution Equations.
 F. Mahmoudi is grateful to SISSA for there kind hospitality.

\end{document}